\documentclass[11pt,reqno]{amsart}
\usepackage[foot]{amsaddr}
\usepackage[top=1in, bottom=1in, left=1.25in, right=1.25in]{geometry}
\usepackage[normalem]{ulem}

\usepackage{amsmath,amsthm,amssymb}
\usepackage{mathrsfs}
\usepackage{mathtools}
\usepackage{thmtools} 
\usepackage{proof} 
\usepackage{enumitem}
\usepackage[dvipsnames]{xcolor} 
\usepackage{graphicx}
\usepackage{booktabs} % makes prettier tables.
\usepackage[colorlinks=true,allcolors=blue,allbordercolors=white]{hyperref} % hyperlink setup
\usepackage{float}

\usepackage{sidecap} % captions on the side
\usepackage{caption}  % caption formatting package
\counterwithin{figure}{section}

\usepackage[toc,page]{appendix}
\usepackage[backend=biber]{biblatex}
\usepackage{biblatex}
\renewcommand\labelenumi{\upshape(\roman{enumi})}
\renewcommand\theenumi\labelenumi

\theoremstyle{definition}
\newtheorem{question}{Question}[section]

\begin{document}

\title{Early Career Researcher Responses to AI in Math}
\author{Anna Brandenberger \and Hanna Mularczyk \and Dora Woodruff}
\address{Department of Mathematics, Massachusetts Institute of Technology}
\email{\{abrande,hannamul,dorawood\}@mit.edu}

\begin{abstract}
We've heard a lot about AI in math from Fields medalists, AI labs, and established professors. However, the personal stakes are highest for early career researchers. Hoping to provide a platform for graduate students and postdocs to voice their opinions, concerns, and experiences, we wrote and distributed a survey on AI in math. 

In this article, we share the responses to our survey. All responses are anonymous, and respondents came from almost $60$ institutions and at least $5$ countries. Our survey included both multiple choice questions and more open-ended free response questions. 

If you are a senior mathematician or university administrator reading these survey results, our goal is for you to understand the concerns of the younger generations and use your authority to take action accordingly. If you are an early career mathematician reading these survey results, we hope you find comfort in the knowledge that your experiences are shared by many others. If you are not a mathematician, we hope you come away from this article with a deeper understanding of the goals, values, and concerns of the community. 
\end{abstract}

\maketitle

\section{Introduction}

In September and early October, we surveyed over $260$ early career mathematicians representing nearly $60$ institutions to get their thoughts on AI in math. Here, we share the survey results. 

\subsection{Motivation}

As AI companies turn to the field of math as a new benchmarking standard, frontier models are getting better at proving mathematical theorems at faster rate than most of us were prepared for:

\begin{enumerate}
    \item In May 2026, an LLM disproved Erdös's unit distance problem \cite{openaierdos}.
    \item In August, ten AI-generated advances in math and computer science were announced
    \cite{openaiten}.
    \item In September, an AI-generated solution to Navier Stokes was announced (as we were preparing the initial distribution of this survey).
    \item In October, right \textit{after} this survey closed to responses, OpenAI released $700+$ solutions to other open problems in math and theoretical computer science.
\end{enumerate} 

As a result, mathematicians are grappling with the fact that a core part of our jobs could soon be automated (\cite{leidendeclaration}, \cite{taoslides}, \cite{weinreichessay}). 

However, the leading voices on the matter have been Fields medalists and famous mathematicians, who can hypothesize and ideate from the safety of a secure job. Us early-career mathematicians do not have the same luxury, and our unique point of view was not as represented in these discussions. Hence, we decided to ask graduate students and postdocs how they felt directly.

This project started as a survey on AI and math specifically for MIT mathematics graduate students, which we released in September $2026$. The results to that survey are discussed in an {open letter} we wrote to the department (which can be viewed at \url{https://math.mit.edu/math-community/math_ai.html}). 

Upon the conclusion of the first survey, we wanted to reach a larger audience and expand on some of the questions (especially: \textit{how} early-career mathematicians are using AI, and the \textit{broader ethical concerns} early-career mathematicians worry about). So, we created and released the second survey, the results of which we share in this article.

\subsection{Survey Methodology}

The survey was administered as a Google form, which opened on September $17$, 2026 and closed on October $8$, 2026. We distributed the survey by asking colleagues and peers at other universities to distribute it to their departments, often through an email list server, on our behalf. The survey also spread through word of mouth and conference and workshop group chats. Anyone who was a graduate student or postdoctoral researcher in mathematics or an adjacent field was encouraged to participate.

\subsection{Disclaimers and caveats}

First and foremost: it is important to note that almost all responses were collected \textit{before} OpenAI released their $700+$ AI-generated problems on October $6$. The responses we collected were not influenced by that event, which may have changed many people's viewpoints. Mathematicians' opinions on the issue are not at all static; the responses in this survey represent only a small window of time. 

This survey \emph{should not be seen as a rigorous scientific data collection effort}. Rather, it is an effort to give early career researchers an organized, anonymous platform to share their experiences and concerns, and to get an informal temperature check. Since each of us has our own experiences and opinions with math and AI, our own biases might have affected the survey results, from the question design to the groups we contacted. 

It should be stated upfront that junior mathematicians 
have very diverse opinions on this issue. While it is important to notice trends and averages in the survey data, it is equally important to not overgeneralize and make imprecise blanket statements about ``what mathematicians believe.'' 
We've made the decision to process the raw survey data very minimally, to avoid filtering it through our own biases. Beyond sorting questions and free responses into broad topic categories, we have not provided summaries or analyses. 

We encourage readers to reach out to the authors with any questions, thoughts, complaints, hesitations, or more. 

\subsection{Article structure}
Our survey contained both multiple choice questions and free response questions. We first share all the data from the multiple choice questions, and then share around \textit{half} of the free responses (selected at random). There were so many free responses that including all of them would make this document unreasonably long. 

In Section~\ref{sec:who}, we give slightly finer information about who responded to our survey. 
In Section~\ref{sec:career}, we share data from multiple choice questions about {job market concerns}. 
In Section~\ref{sec:how-use}, we share data about how respondents use AI in their day-to-day activities. 
In Section~\ref{sec:research-impact}, we share data about how AI has impacted respondents' research activities. 
In Section~\ref{sec:longterm}, we share data concerning how respondents feel about the longterm health of mathematics. We personally found the responses to this section to be the most interesting. 
In Section~\ref{sec:concluding-q}, we share data from concluding questions (for example: how much do researchers think about ethical issues beyond math when they consider using AI?). 
Finally, in the appendix, we share a selection of the free responses, broken up into $6$ subsections by question. 

\subsection*{Acknowledgments} We want to give a personal word of thanks to our respondents.
Being a graduate student or postdoc in math was never an easy job. We're expected to do a lot: research, teach, mentor, learn, present work... Now, on top of all that, we're being forced to reckon with existential questions about a field we love, all while incredibly uncertain about our own futures. Whether you are optimistic, pessimistic, or just confused about AI in math, I think we can all agree: the past couple of months have been \textit{exhausting}. 

With that in mind, we want to thank the $268$ graduate students and postdocs who took the time and energy to fill out our survey. Many of them left long, thoughtful paragraphs in the free response section and brought up points we had not considered. 
Many others volunteered to help distribute the survey within their departments. Without their help, we could not have reached so many people. 

\section{Who Responded to Our Survey?}\label{sec:who}

In total, $268$ mathematicians responded to our survey. Of these, around $85$ percent were graduate students and around $15$ percent were postdocs. The majority identified pure math as their primary area, but we also had responses from applied mathematicians, theoretical computer scientists, and statisticians: 
\begin{figure}[H]
    \centering
    \vspace{-.5\baselineskip}
    \includegraphics[width=.31\linewidth]{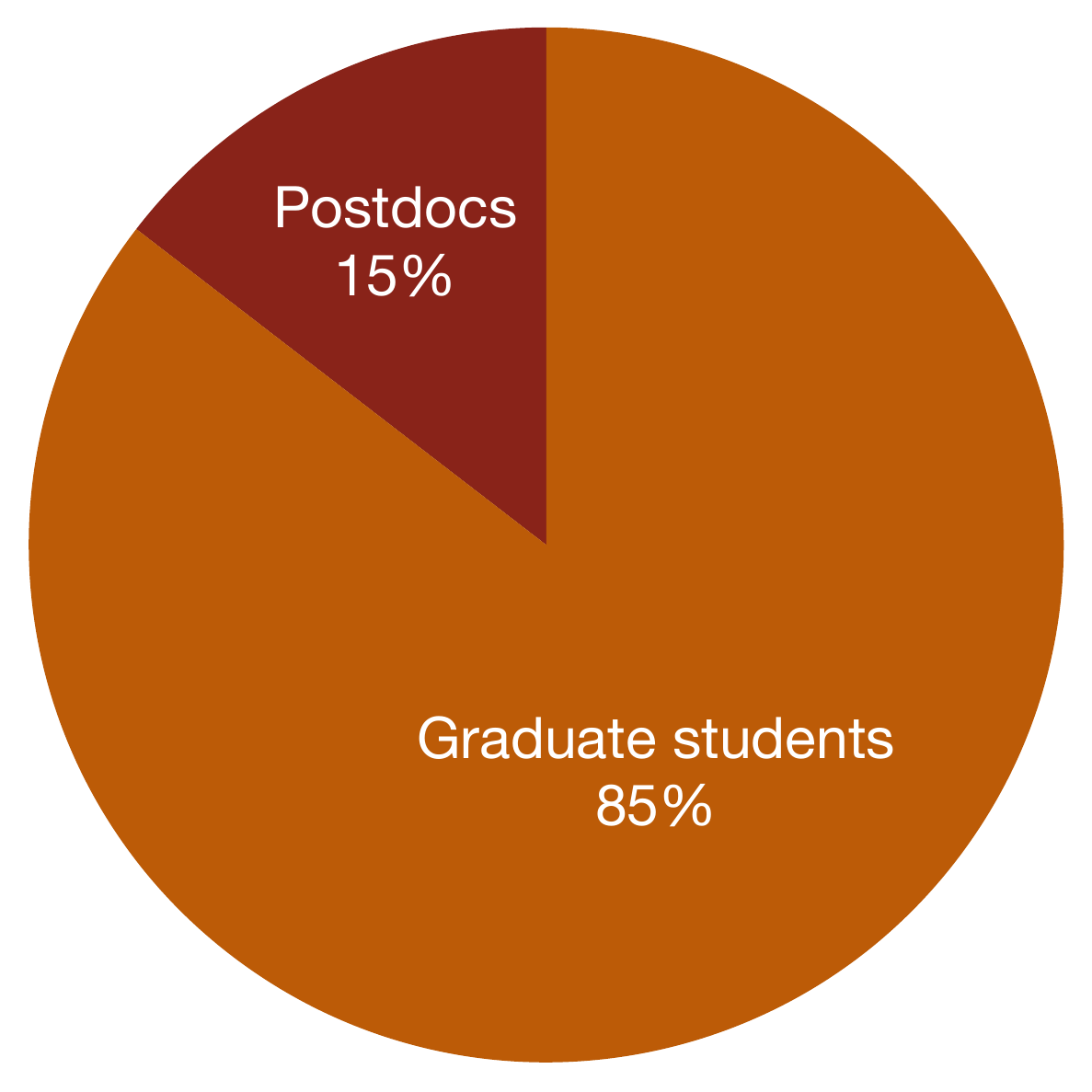}
    \hspace{.4cm}
    \includegraphics[width=.55\linewidth]{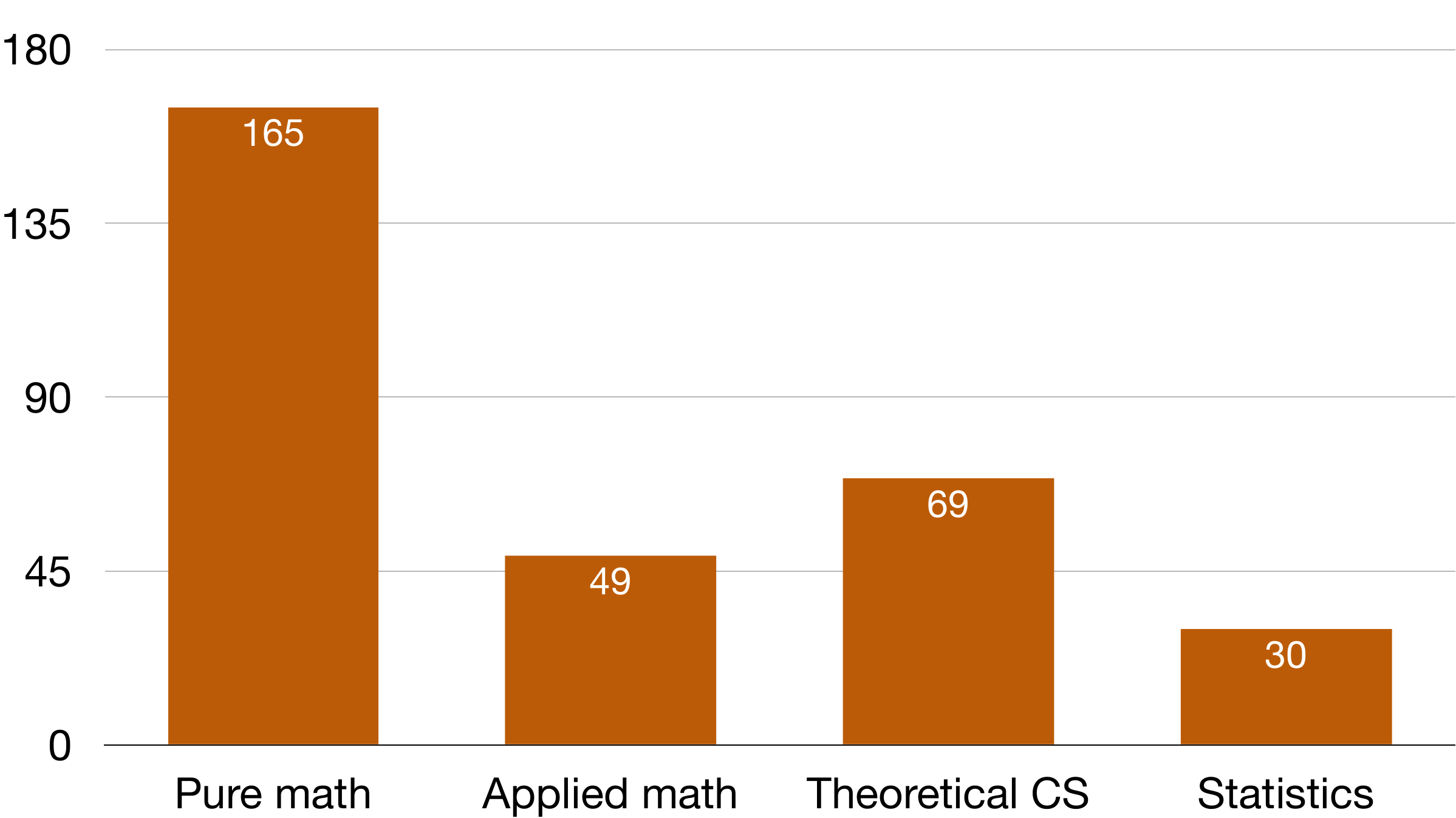}
    \caption{Breakdown of respondents by position and field (we let respondents choose more than one field if applicable). ($n=268$)}
    \label{fig:status}
\end{figure}

Respondents came from nearly $60$ different institutions, including both public and private universities. The majority of respondents were based in the US, but we also received responses from Canada, India, the Netherlands, Germany, England, and France. 
As several respondents shared concerns that their affiliations would be identifying, we won't make the full list of represented institutions public. Please contact the authors if you have further questions about who participated in the survey.

\section{Career and Job Market Concerns}
\label{sec:career}
\begin{question}\label{q:future-plans}
    Have you considered changing your future career/academic plans on account of AI?
\end{question}

\begin{figure}[H]
    \centering
    \vspace{-1.5\baselineskip}
    \includegraphics[width=.45\linewidth]{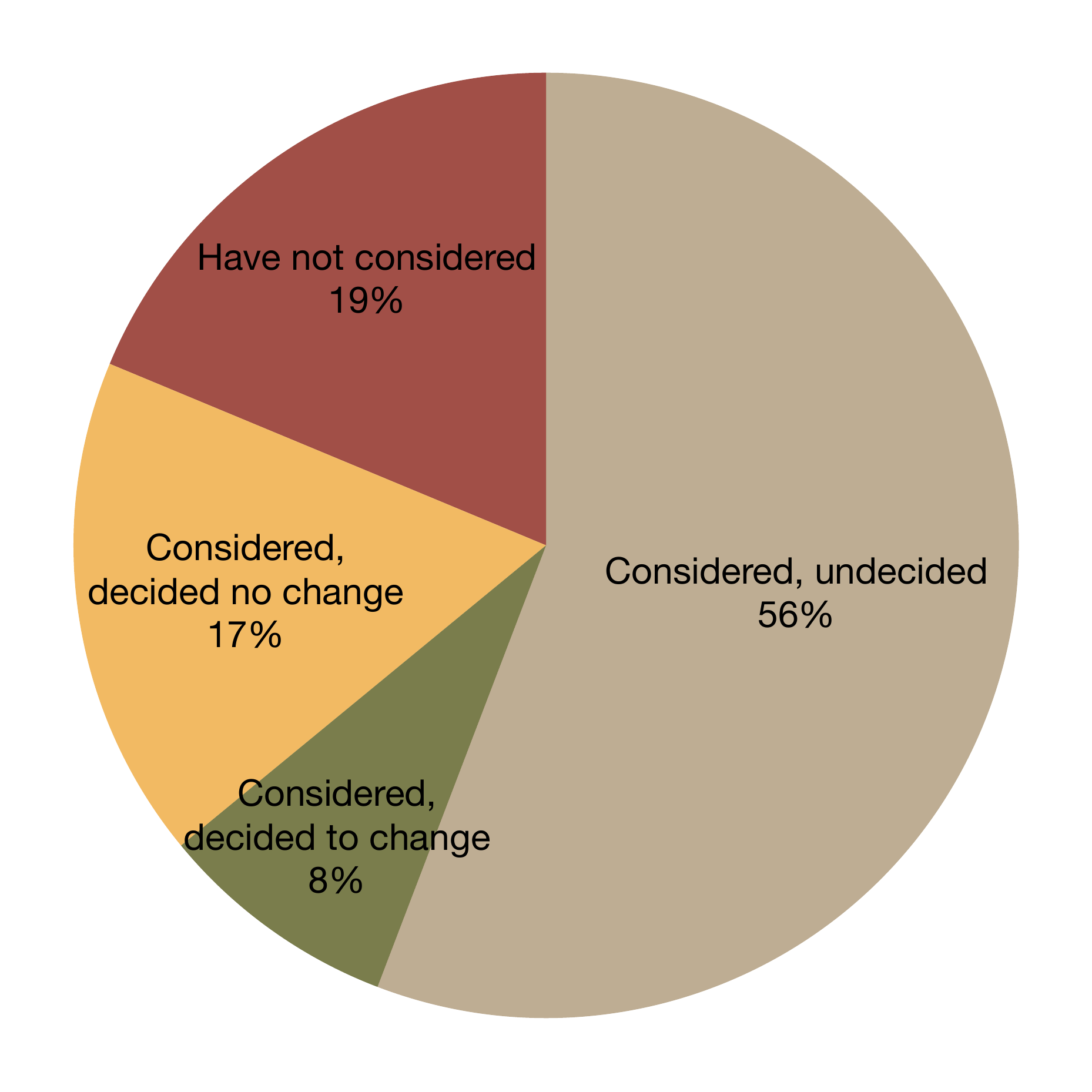}
    \vspace{-1.5\baselineskip}

    \caption{Responses to Question~\ref{q:future-plans} ($n=267$).}
    \label{fig:career-change}
\end{figure}

\begin{question}\label{q:continue}
     Concretely, has the probability of you continuing in academic math after your PhD/postdoc increased or decreased over the last couple of months? 
\end{question}

\begin{figure}[H]
    \centering    
    \vspace{-\baselineskip}
    \includegraphics[width=\linewidth]{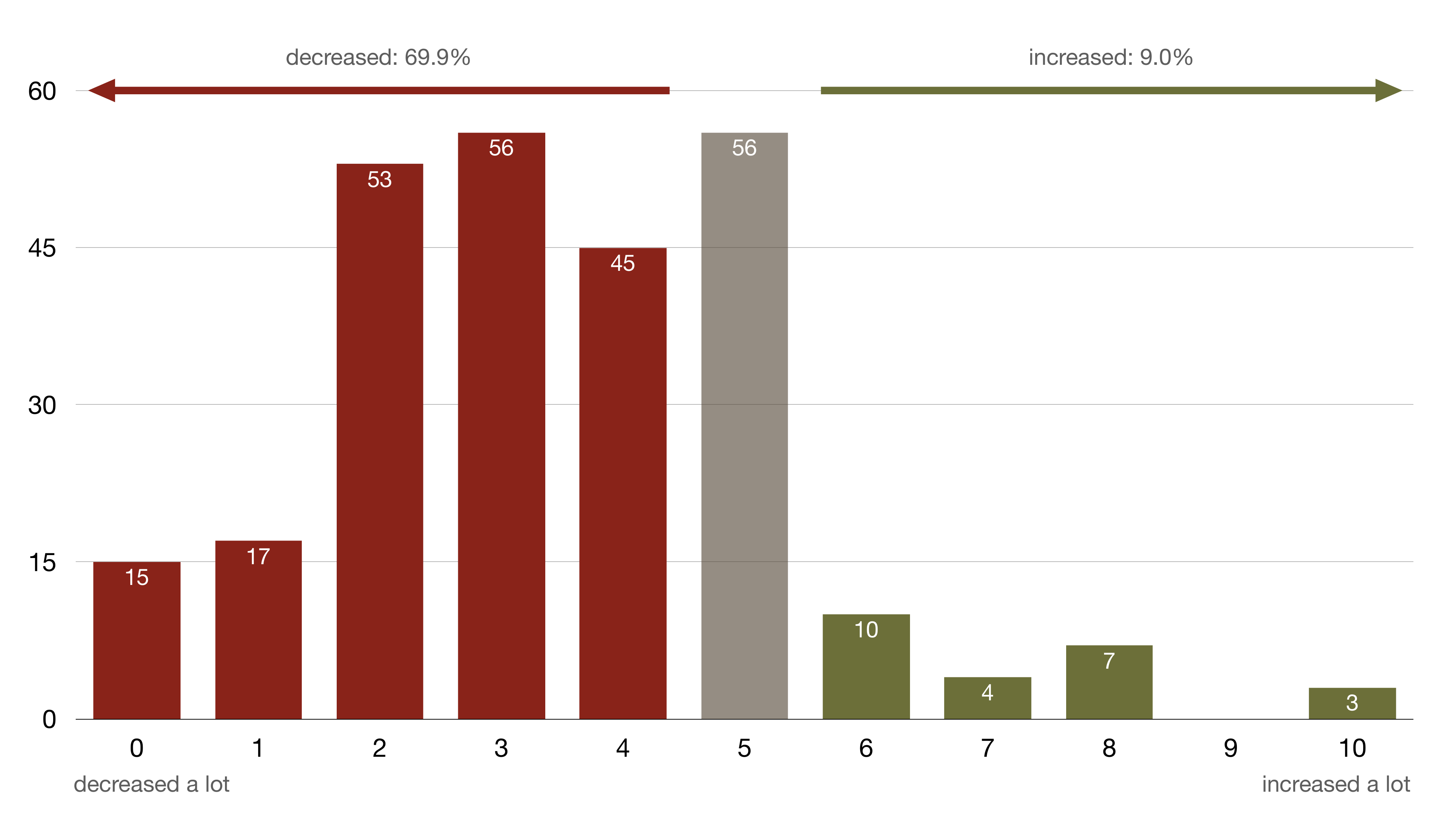}
    \vspace{-2\baselineskip}
    \caption{Responses to Question~\ref{q:continue} ($n=266$).}
    \label{fig:continue}
\end{figure}

\begin{question}\label{q:tenure-track}
    Do you think AI developments will make it harder or easier for you to get a tenure-track job?
\end{question}

\begin{figure}[H]
    \centering
    \vspace{-\baselineskip}
    \includegraphics[width=\linewidth]{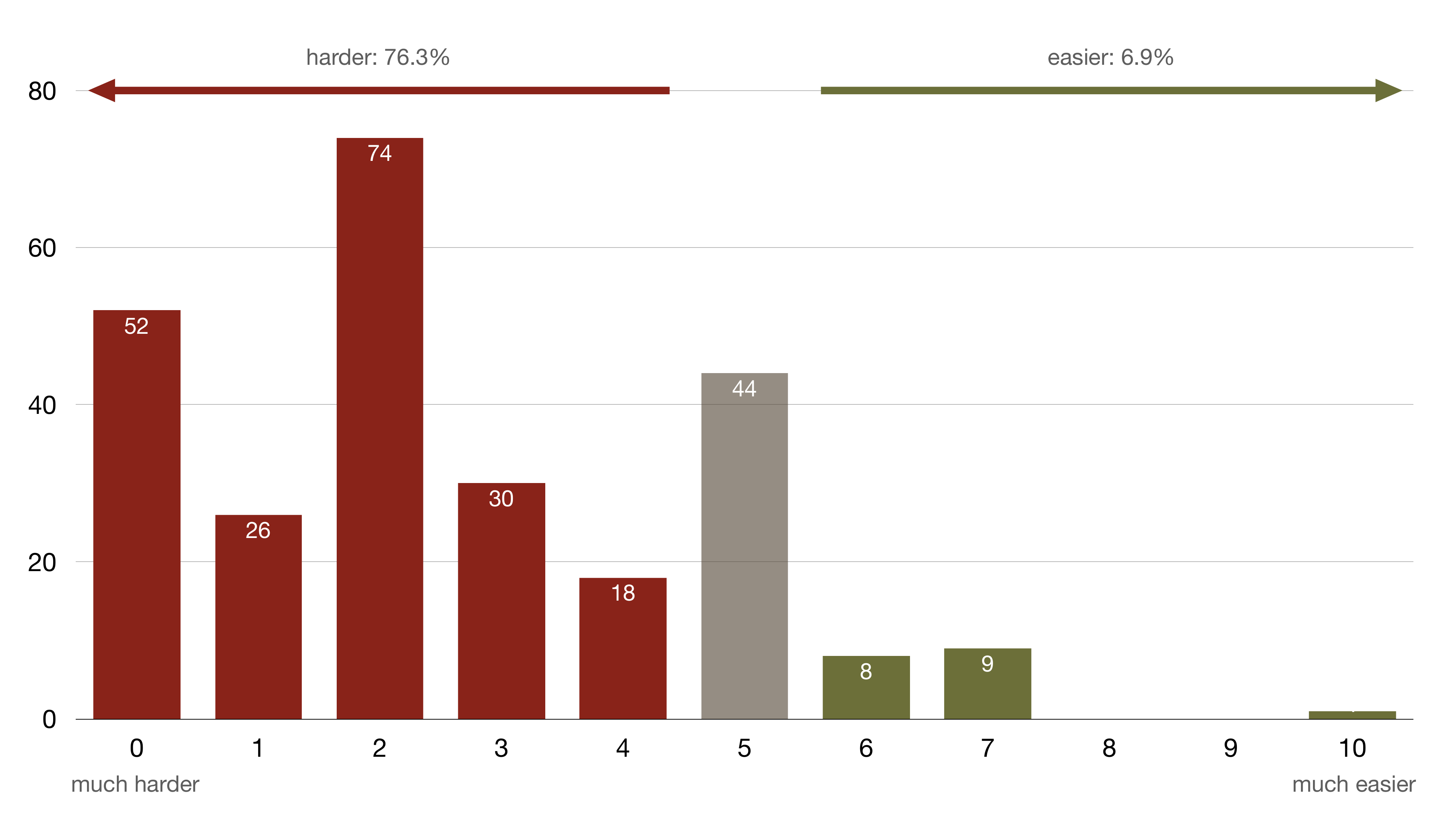}
    \vspace{-2\baselineskip}
    \caption{Responses to Question~\ref{q:tenure-track} ($n=262$).}
    \label{fig:tenure-track}
\end{figure}

\clearpage

\section{How do respondents use AI?}
\label{sec:how-use}
\begin{question}\label{q:howmuch}
    How much do you use AI for your math research?
\end{question}

\begin{figure}[H]
    \centering
    \vspace{-\baselineskip}
    \includegraphics[width=\linewidth]{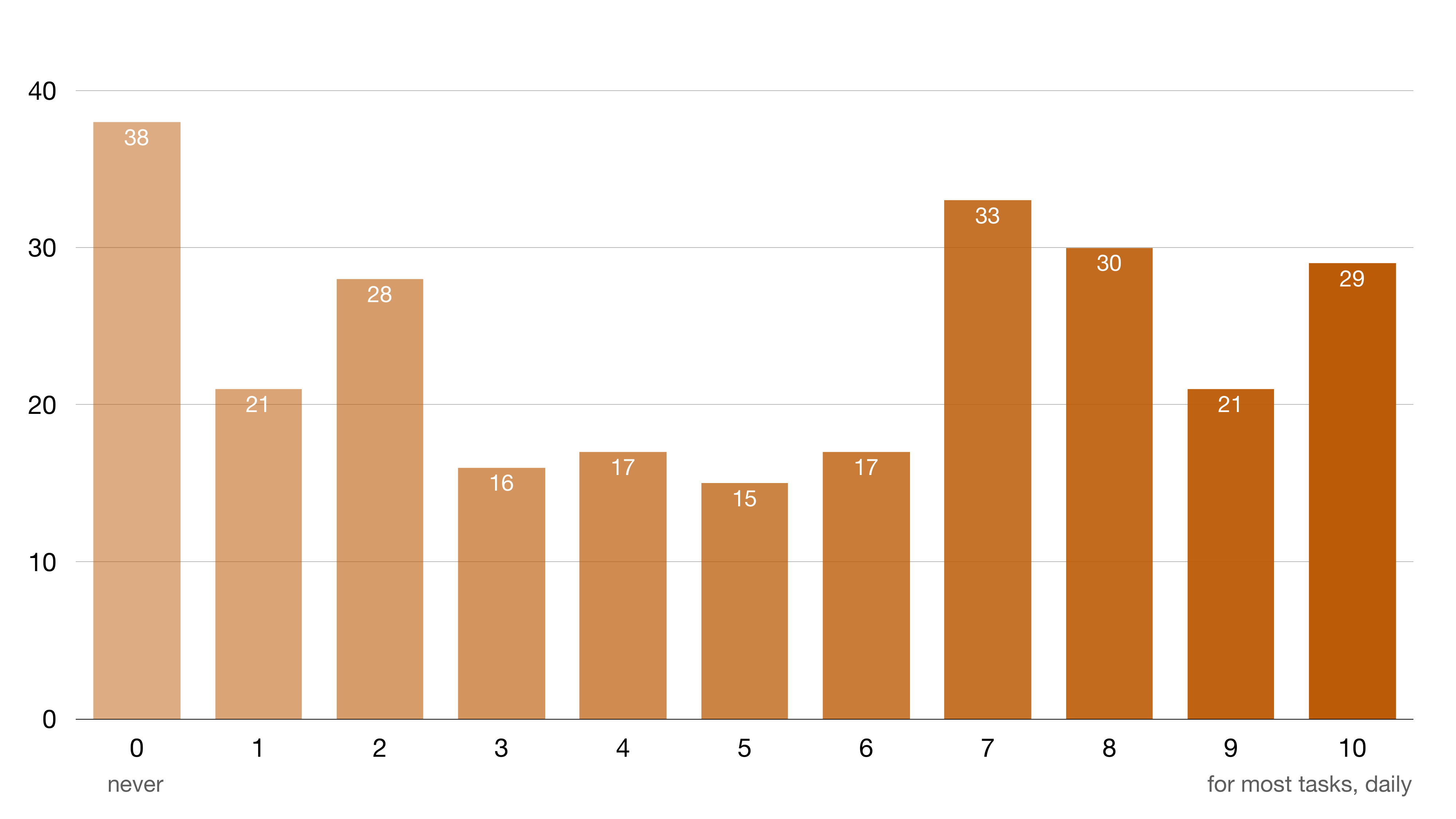}
    \vspace{-2\baselineskip}
    \caption{Responses to Question~\ref{q:howmuch} ($n=265$).}
    \label{fig:ai-use}
\end{figure}

\begin{question}\label{q:enjoy-ai}
    If you use AI for your research, do you enjoy doing so?
\end{question}

\begin{figure}[H]
    \centering
    \vspace{-.5\baselineskip}
    \includegraphics[width=\linewidth]{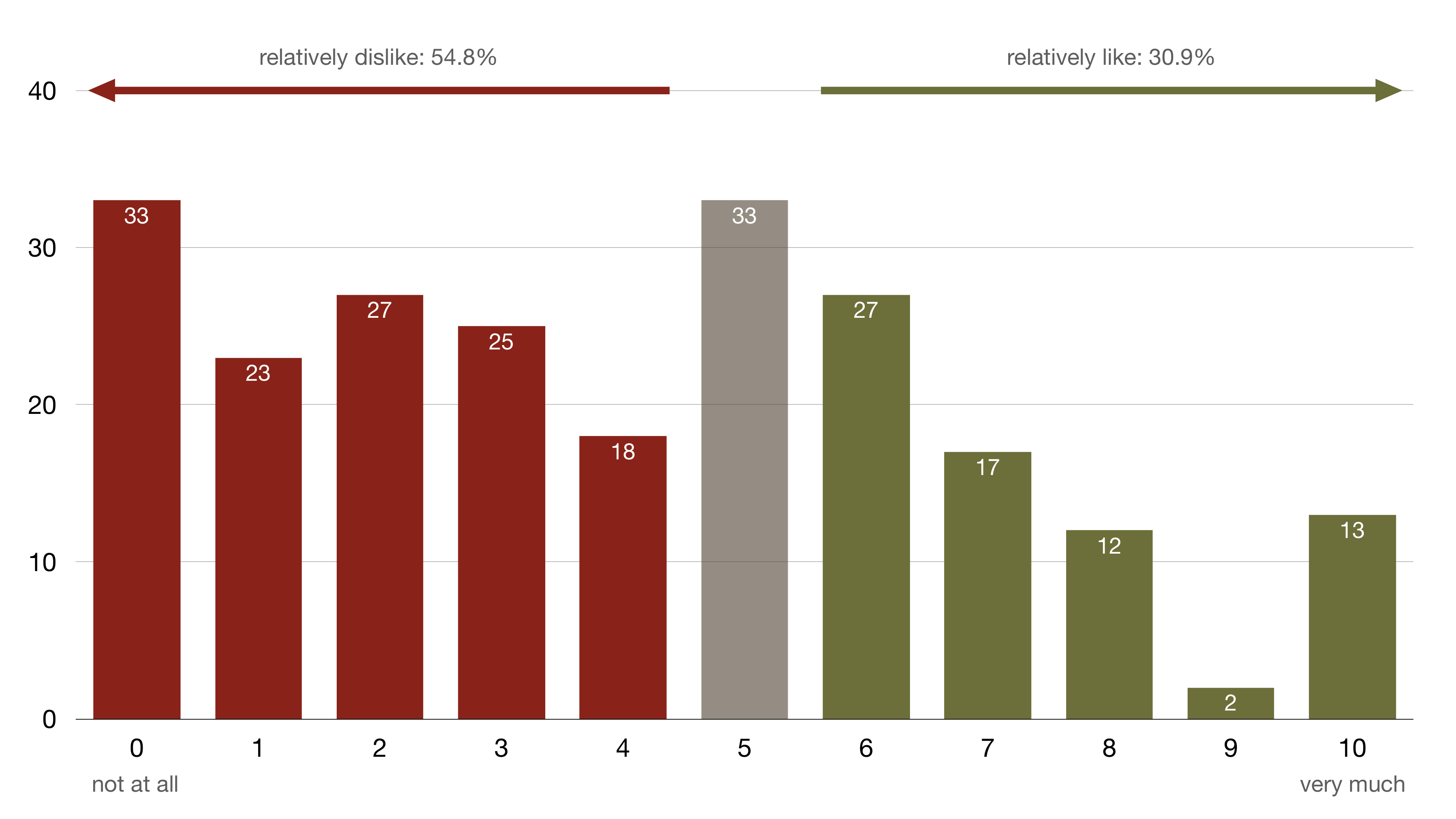}
    \vspace{-2\baselineskip}
    \caption{Responses to Question~\ref{q:enjoy-ai} ($n=230$).}
    \label{fig:enjoy}
\end{figure}

\begin{question}\label{q:which-tasks}
    What kinds of tasks do you use AI for, if any? (We provided a checklist of options, and the ability for respondents to write in their own response. Respondents were invited to check off all applicable answers, so percentages are out of the number of respondents who answered this question.)
\end{question}

\begin{figure}[H]
    \centering
    \vspace{-\baselineskip}
    \includegraphics[width=\linewidth]{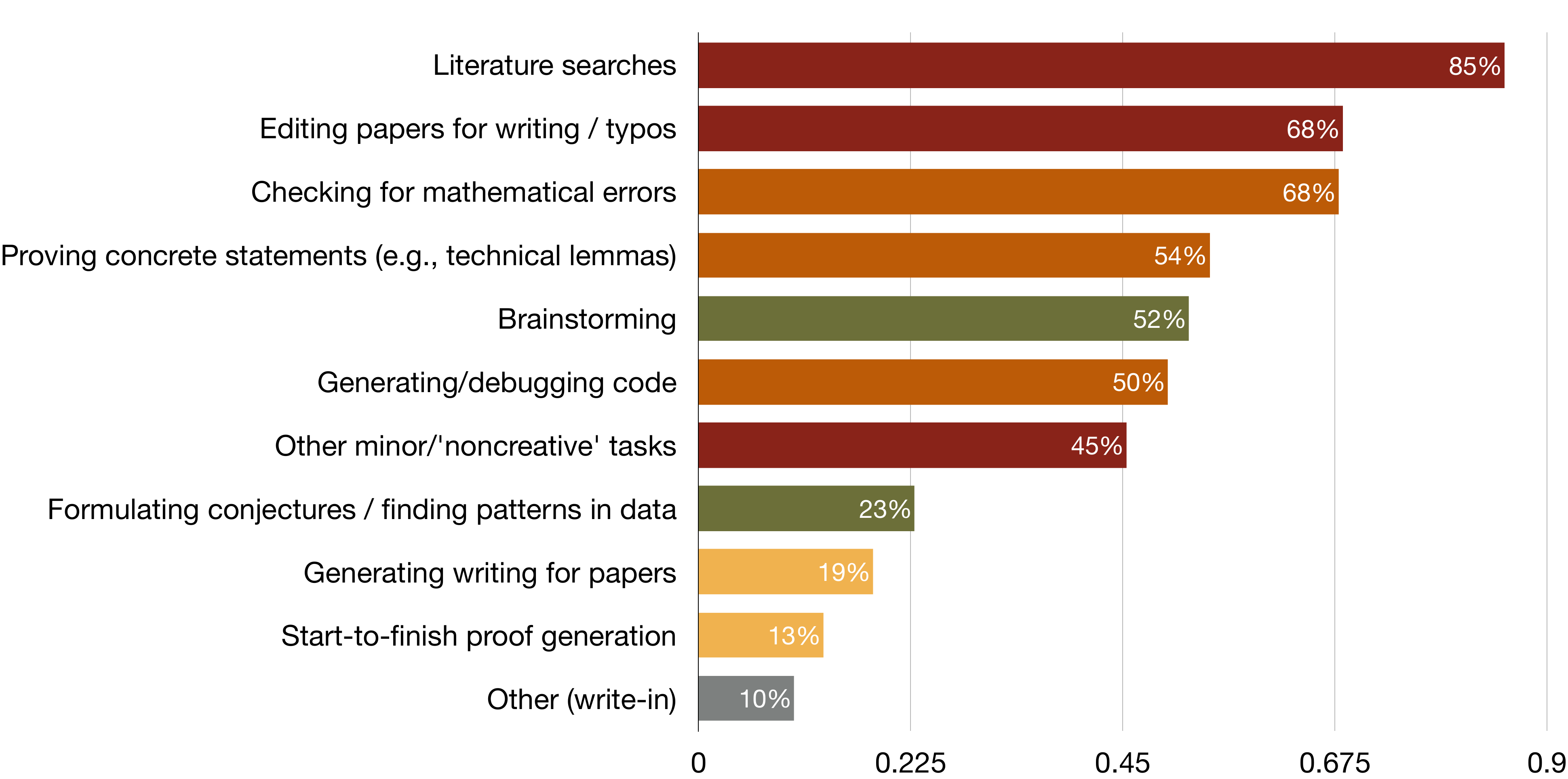}
    \caption{Responses to Question~\ref{q:which-tasks}. 
    Among the Other (write-in) responses, several (13) described various ways they use AI to help explore and understand literature. (Percentages out of $n=227$; the 15\% respondents who did not answer this question mostly use AI minimally).}
    \label{fig:tasks}
\end{figure}

\begin{question}\label{q:which-models}
    Which models, if any, do you use? 
\end{question}
\begin{figure}[H]
    \centering
    \vspace{-\baselineskip}
    \includegraphics[width=\linewidth]{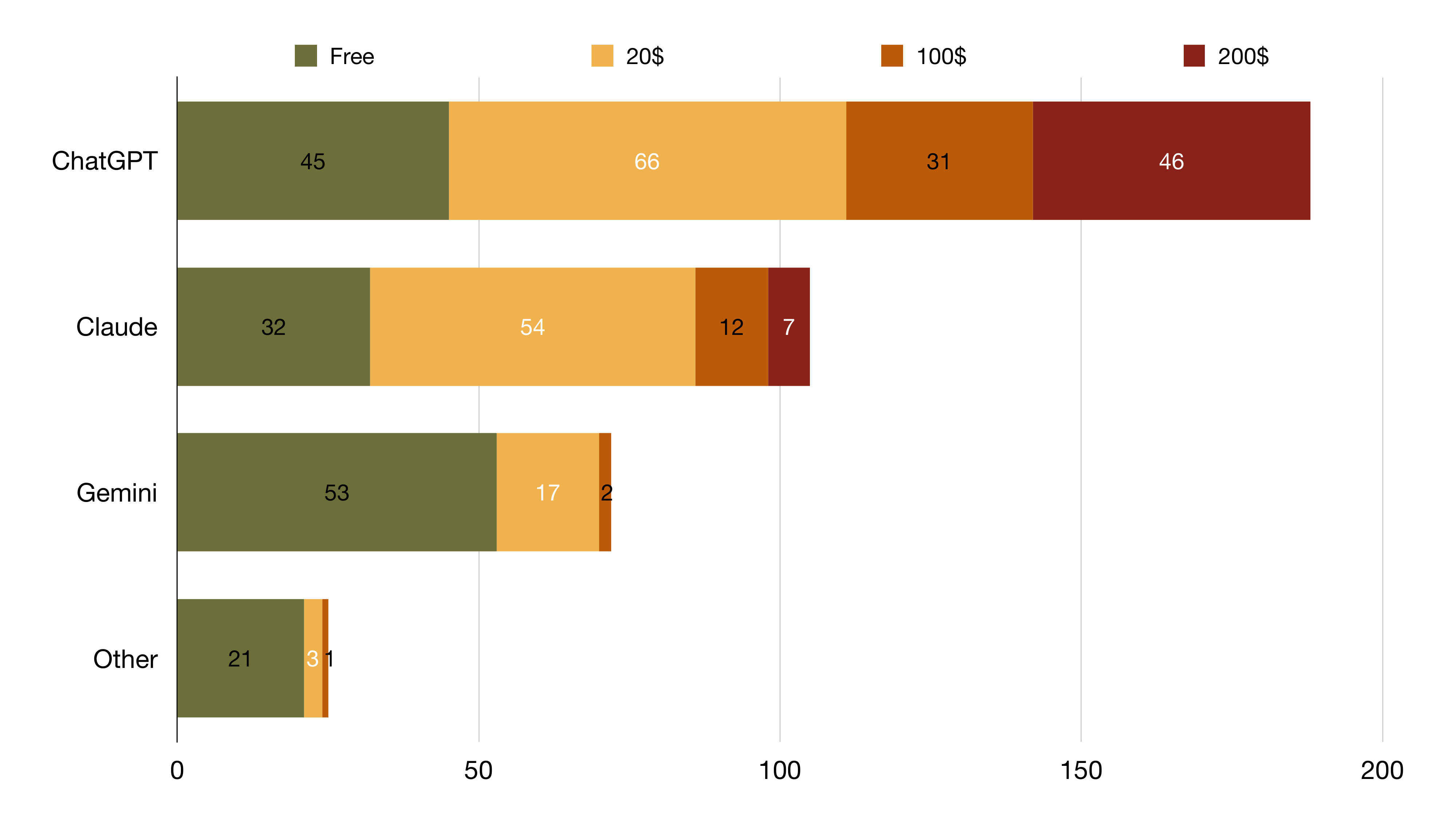}
    \vspace{-2\baselineskip}
    \caption{Responses to Question~\ref{q:which-models} ($n=268$).}
    \label{fig:which-model}
\end{figure}

\begin{question}
    If you use a subscription, how do you pay for it? (Respondents were invited to check off all applicable answers, so percentages are out of $n=181$.)
\end{question}
\begin{enumerate}
    \item Out of pocket: 62\%
    \item Advisor / research group: 20\% 
    \item Free access from AI lab: 12\% 
    \item Department / university / ``sandbox" environment: 23\% 
\end{enumerate}

\section{How has AI impacted respondents' research?}
\label{sec:research-impact}

\begin{question}\label{q:ai-solving}
     How would you feel about AI solving a problem you've worked on? (Either a lab, or an individual one-shotting a problem) 
\end{question}

\begin{figure}[H]
    \centering
    \vspace{-\baselineskip}
    \includegraphics[width=\linewidth]{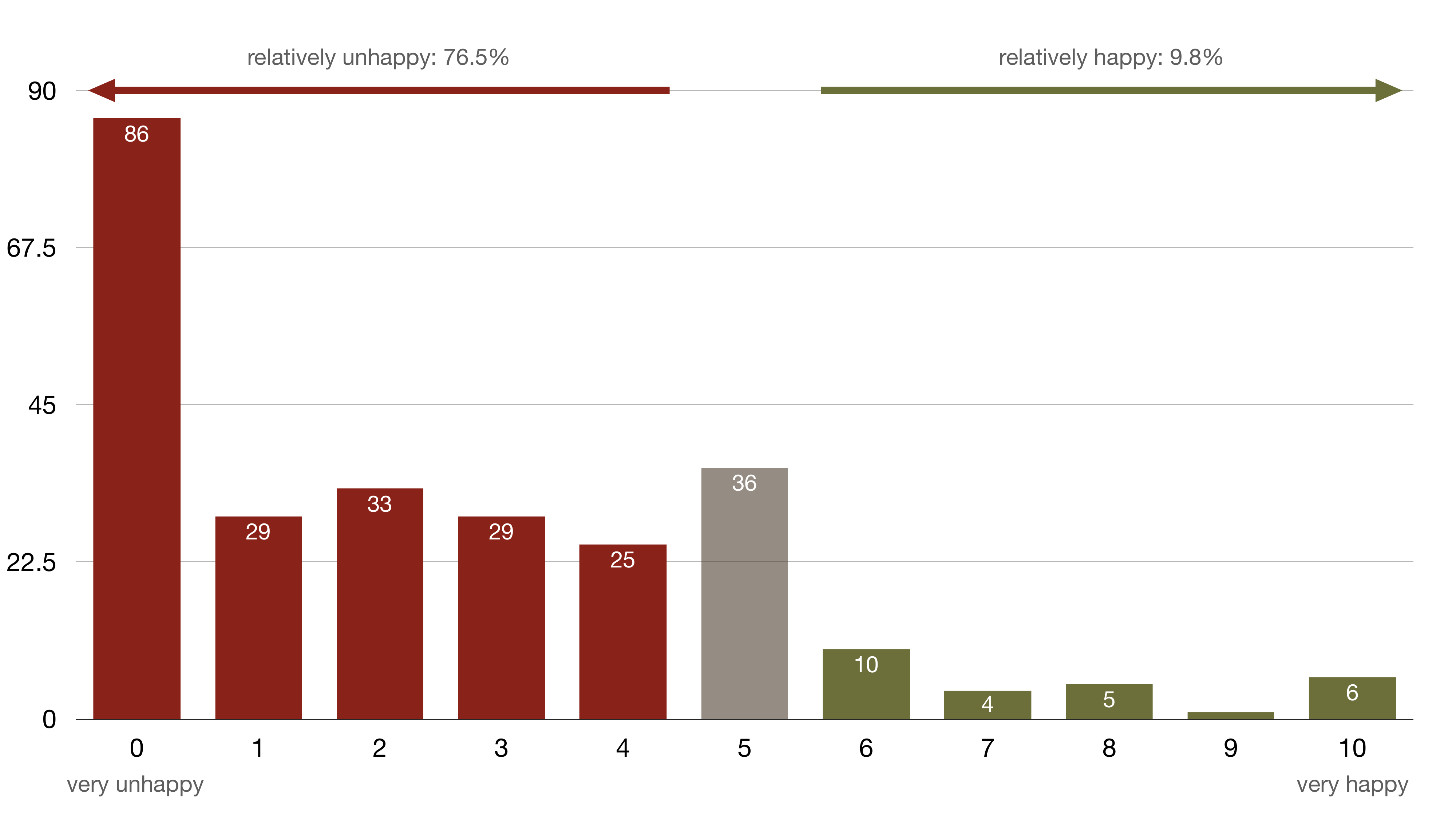}
    \vspace{-2\baselineskip}
    \caption{Responses to Question~\ref{q:ai-solving} ($n=264$).}
    \label{fig:ai-solving}
\end{figure}

\begin{question}\label{q:motivation}
    Do you feel more or less motivated to work on research as a result of AI developments?
\end{question}

\begin{figure}[H]
    \centering
    \vspace{-\baselineskip}
    \includegraphics[width=\linewidth]{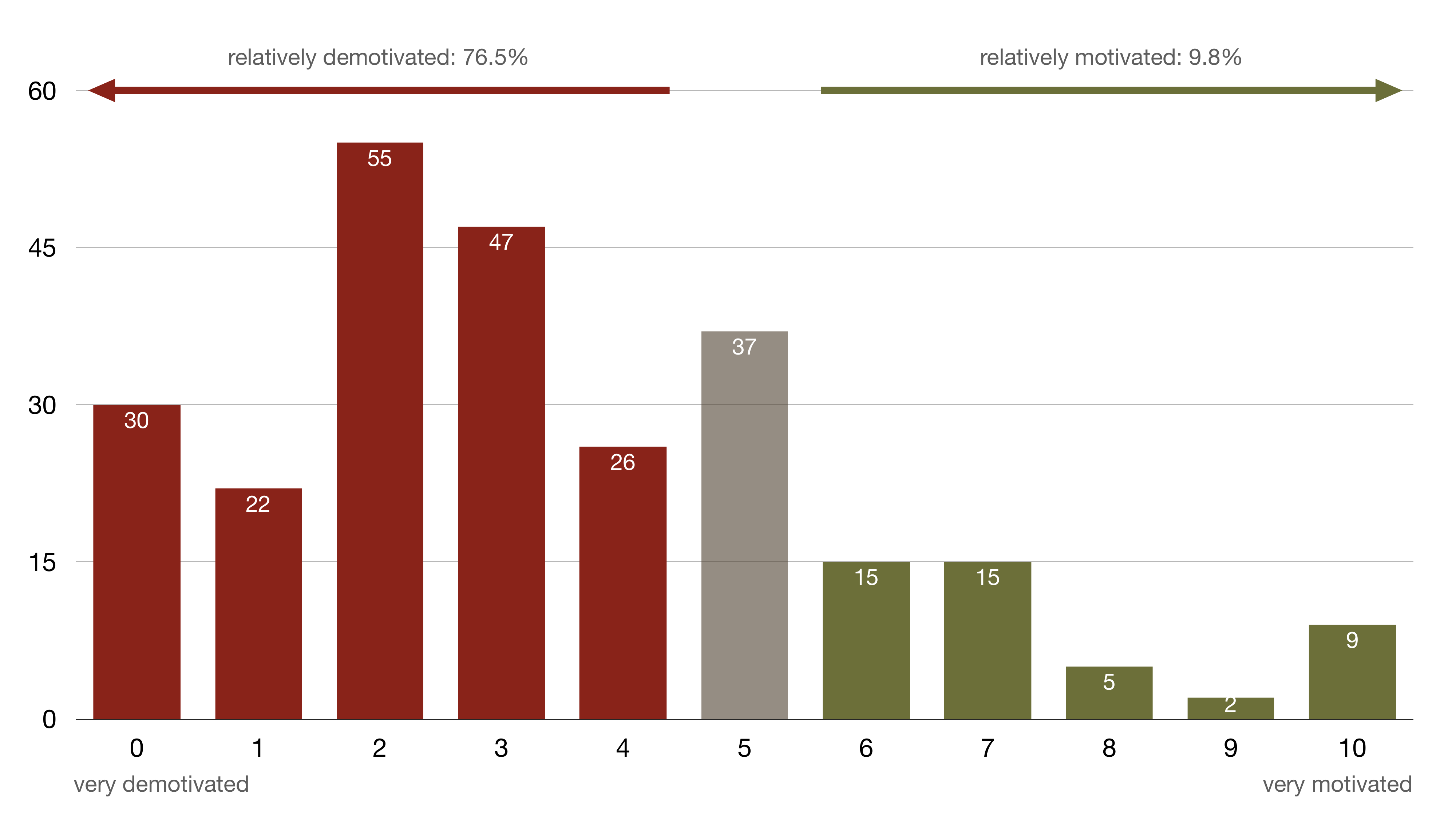}
    \vspace{-2\baselineskip}
    \caption{Responses to Question~\ref{q:motivation}. ($n=263$)}
    \label{fig:motivation}
\end{figure}

\begin{question}\label{q:consequences}
    Has AI already had concrete consequences for your personal growth as a mathematician? Would you consider those effects positive or negative? 
\end{question}

\begin{figure}[H]
    \centering
    \includegraphics[width=0.45\linewidth]{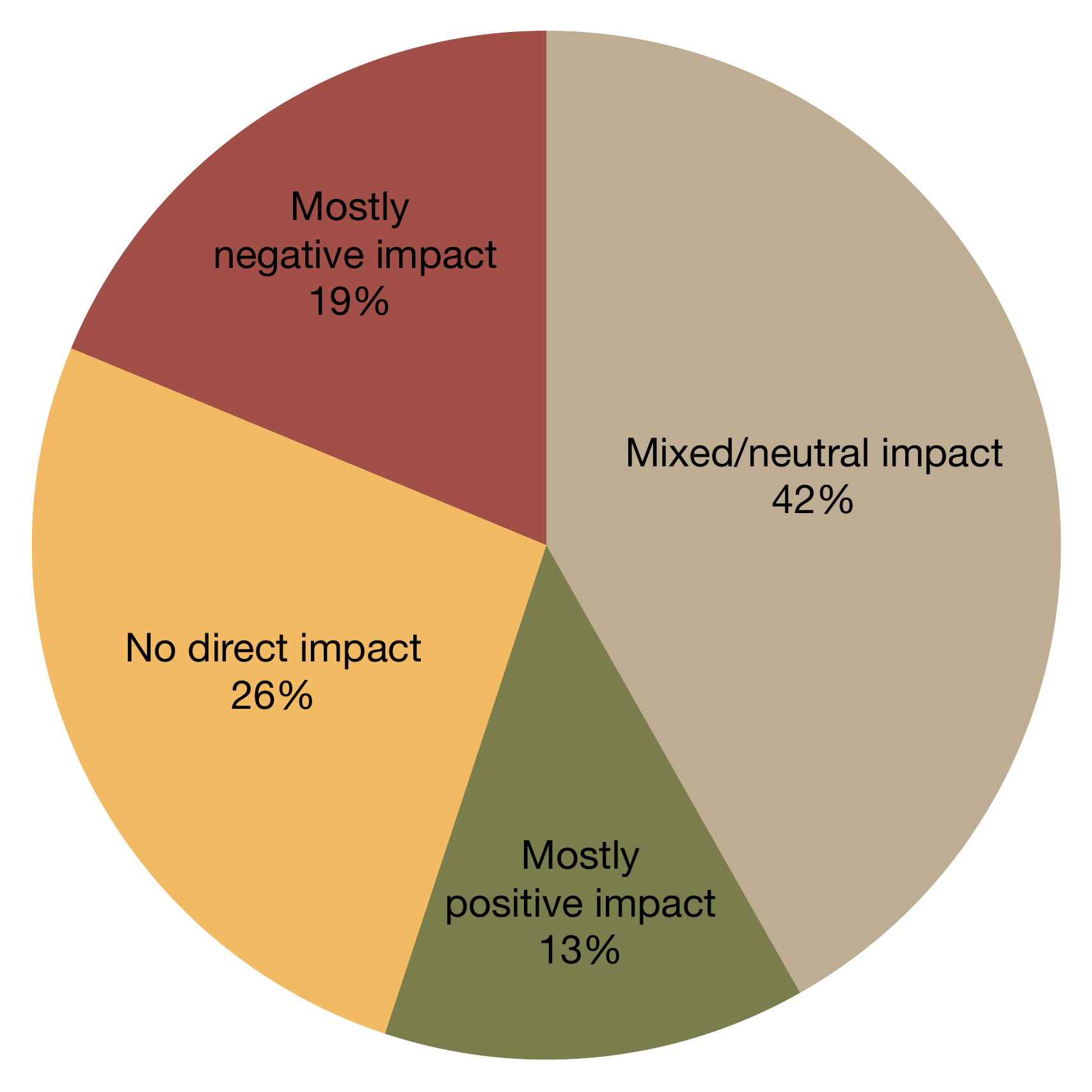}
    \caption{Responses to Question~\ref{q:consequences} ($n=256$)}
    \label{fig:consequences}
\end{figure}

\section{Concerns about the longterm health of the field}
\label{sec:longterm}

\begin{question}\label{q:fave}
     Which would you say are your favourite aspects of being in math? Why do you do this? We know it's hard to prioritize, but please pick up to five options (to avoid everyone clicking everything).
\end{question}

\begin{figure}[H]
    \centering
    \includegraphics[width=\linewidth]{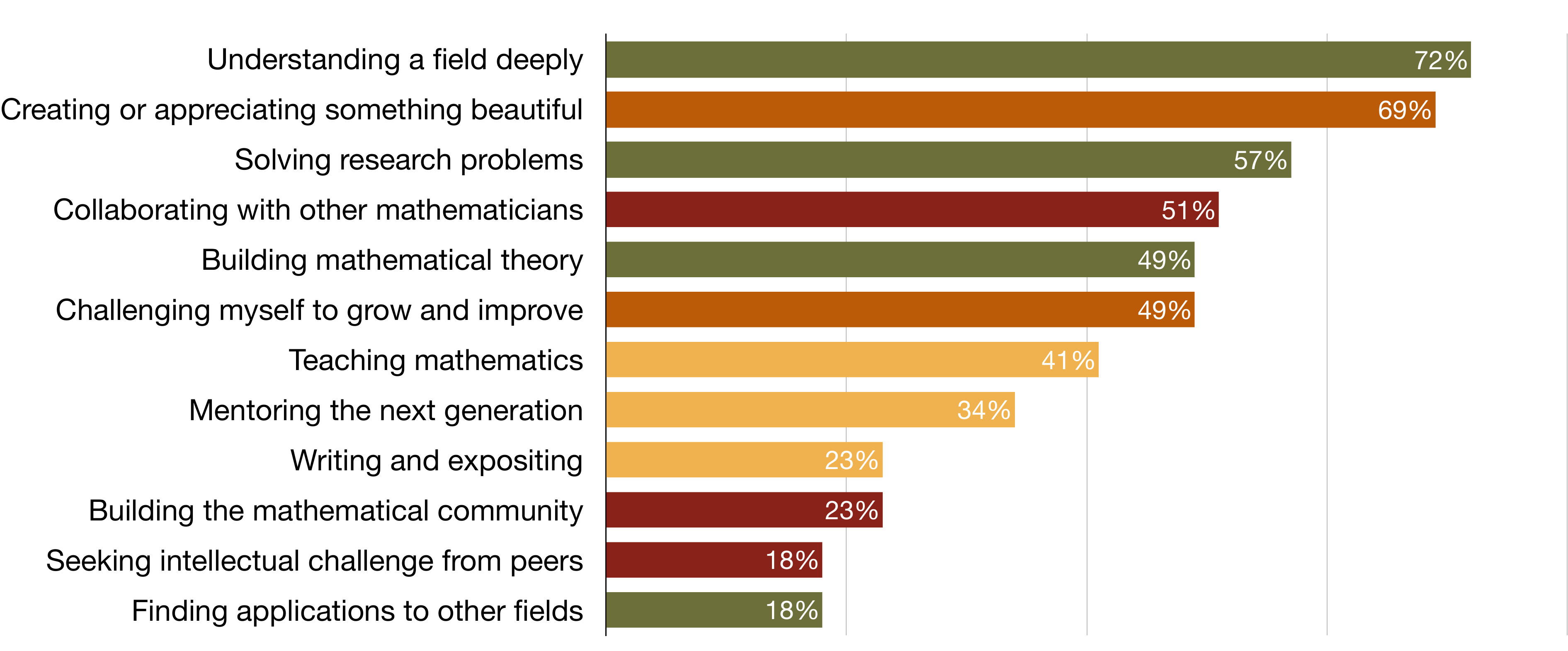}
    \caption{Responses to Question~\ref{q:fave} (Percentages out of $n=265$).}
    \label{fig:fave}
\end{figure}

\begin{question}\label{q:digestion}
    If ``proof digestion" became the main creative task of human mathematicians, would you still want to pursue a career in mathematics? (By proof digestion, we mean `everything that comes between formal validation a proof, and the general community deeply understanding the proof'; e.g., giving talks, writing surveys, conversing with others, contributing to textbooks). 
\end{question}

\begin{figure}[H]
    \centering
    \vspace{-\baselineskip}
    \includegraphics[width=\linewidth]{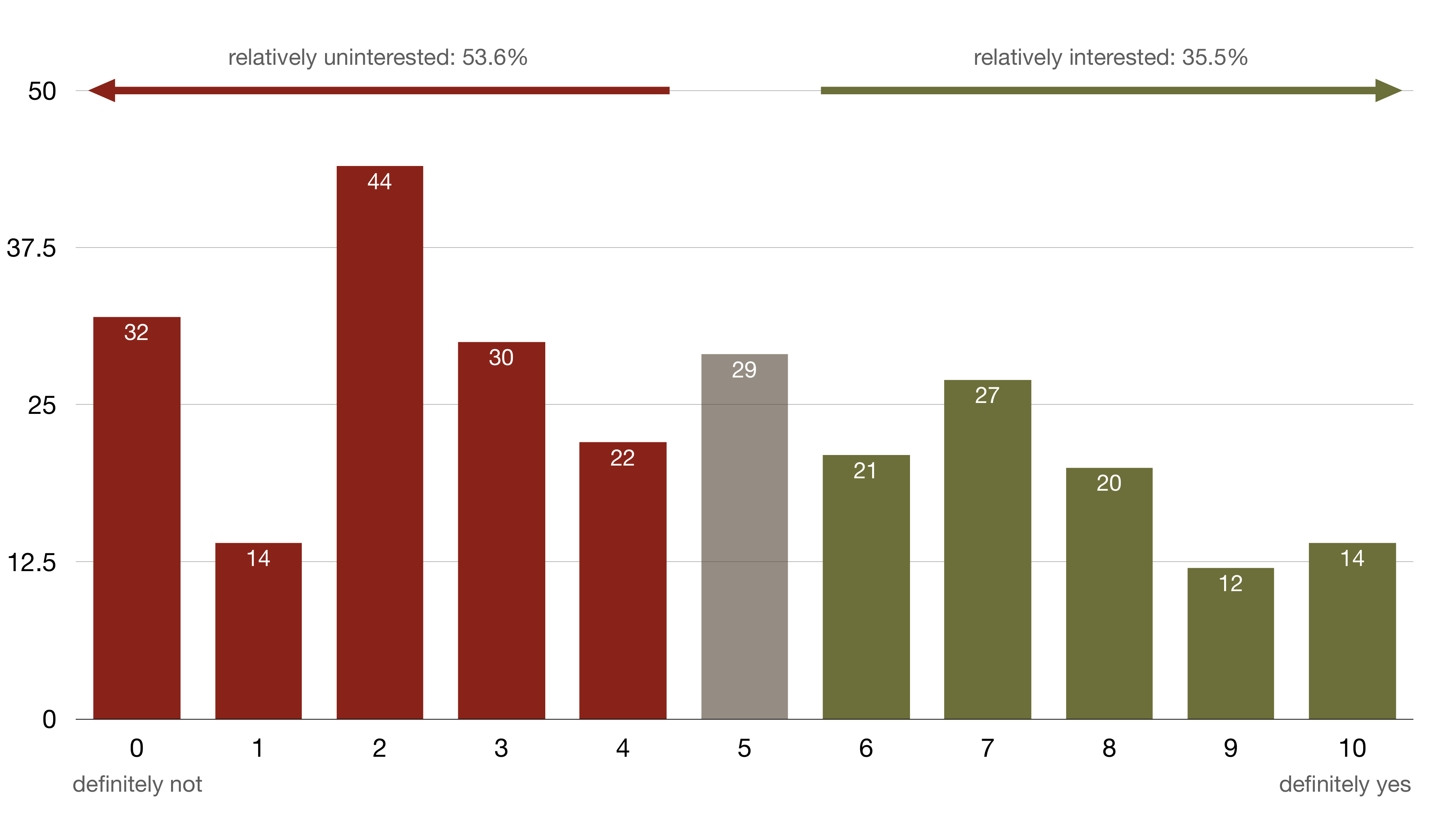}
    \vspace{-2\baselineskip}
    \caption{Responses to Question~\ref{q:digestion} ($n=265$).}
    \label{fig:digestion}
\end{figure}

\begin{question}\label{q:issues}
    Many potential issues have been raised by the community. Which, if any, worry you \textit{significantly}? Respondents were allowed to select many options.
\end{question}

\begin{figure}[H]
    \centering
    \vspace{-\baselineskip}
    \includegraphics[width=\linewidth]{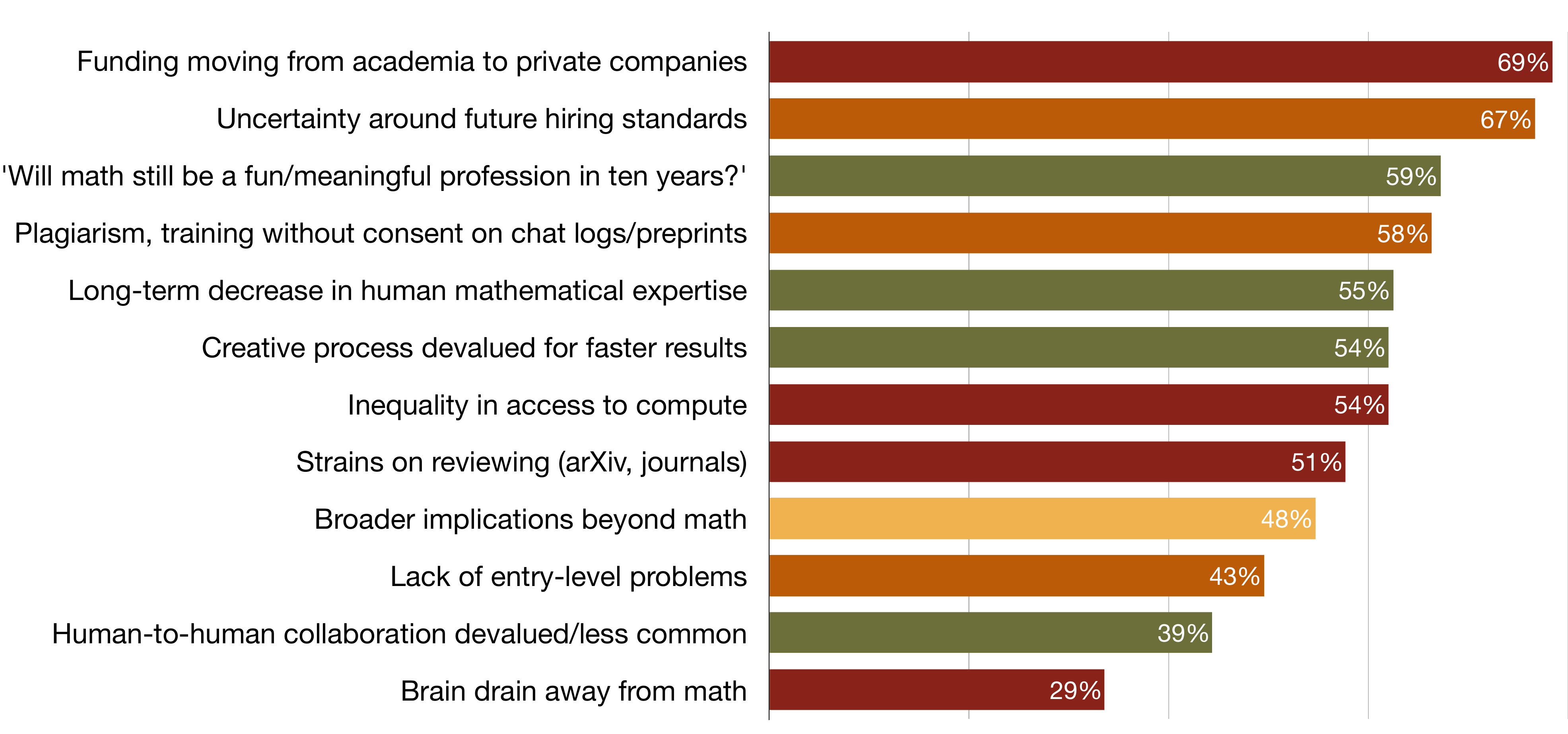}
    \vspace{-1.5\baselineskip}
    \caption{Responses to Question~\ref{q:issues}. (Percentages out of $n=265$.) The number of answers were correlated with respondents' general opinion on AI in math (Question~\ref{q:overall}) below. Broader implications beyond math, strains on reviewing and inequality in access to compute were the concerns that all groups \emph{agreed on the most}. }
    \label{fig:concerns}
\end{figure}

\section{Concluding questions}
\label{sec:concluding-q}

\begin{question}\label{q:ethics}
    To what extent do you keep in mind broader, non-mathematical ethical issues in mind (e.g., misinformation spread, military contracts) when you consider using AI for research?
\end{question}

\begin{figure}[H]
    \centering
    \vspace{-.5\baselineskip}
    \includegraphics[width=\linewidth]{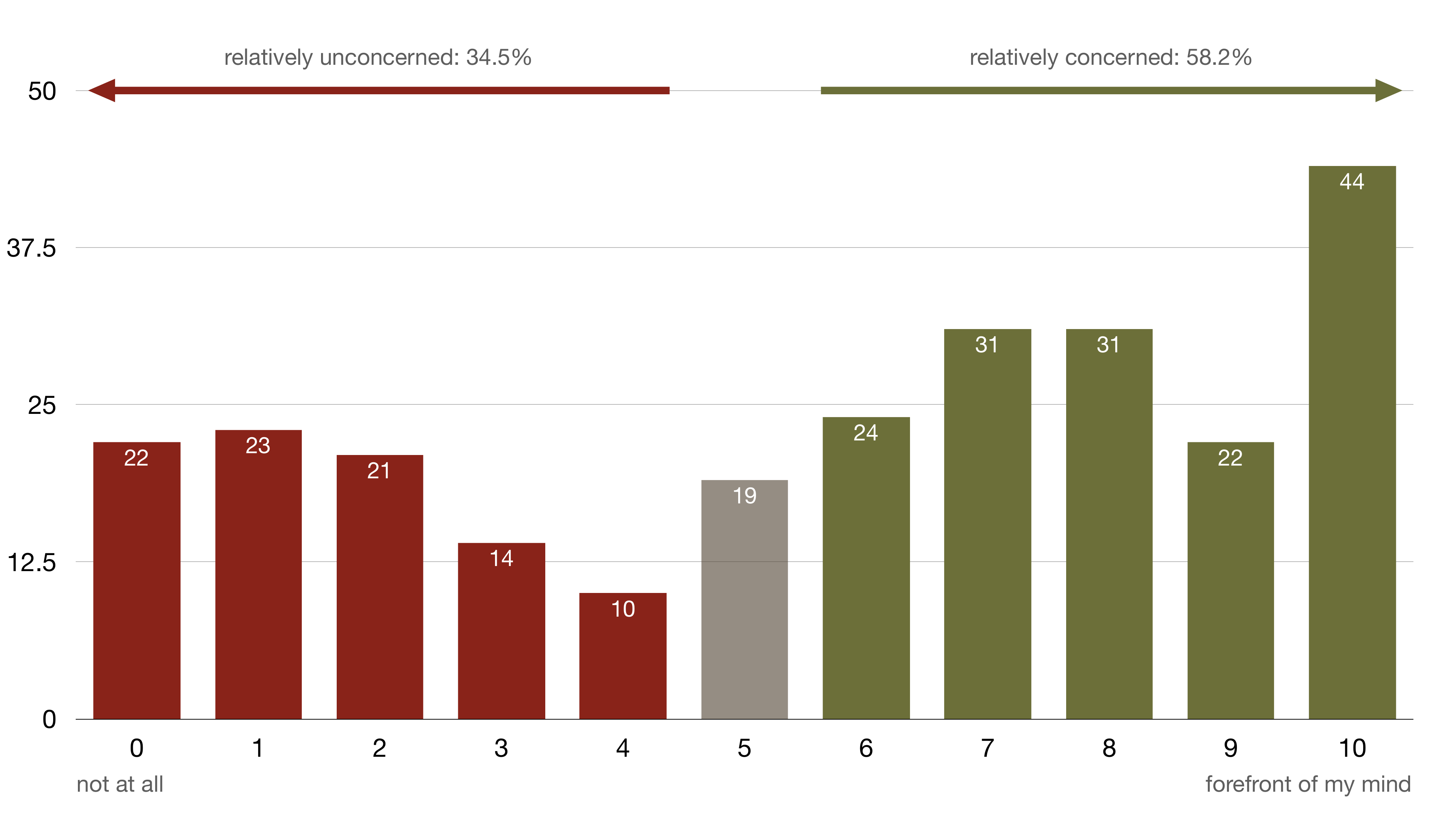}
    \vspace{-2\baselineskip}
    \caption{Responses to Question~\ref{q:ethics} ($n=261$).}
    \label{fig:ethics}
\end{figure}

\begin{question}\label{q:declaration}
    Have you seen the declaration `A Severe Misalignment of AI in Mathematics?’ Do you have an opinion on it?
\end{question}

\begin{figure}[H]
    \centering
    \vspace{-2\baselineskip}
    \includegraphics[width=.55\linewidth]{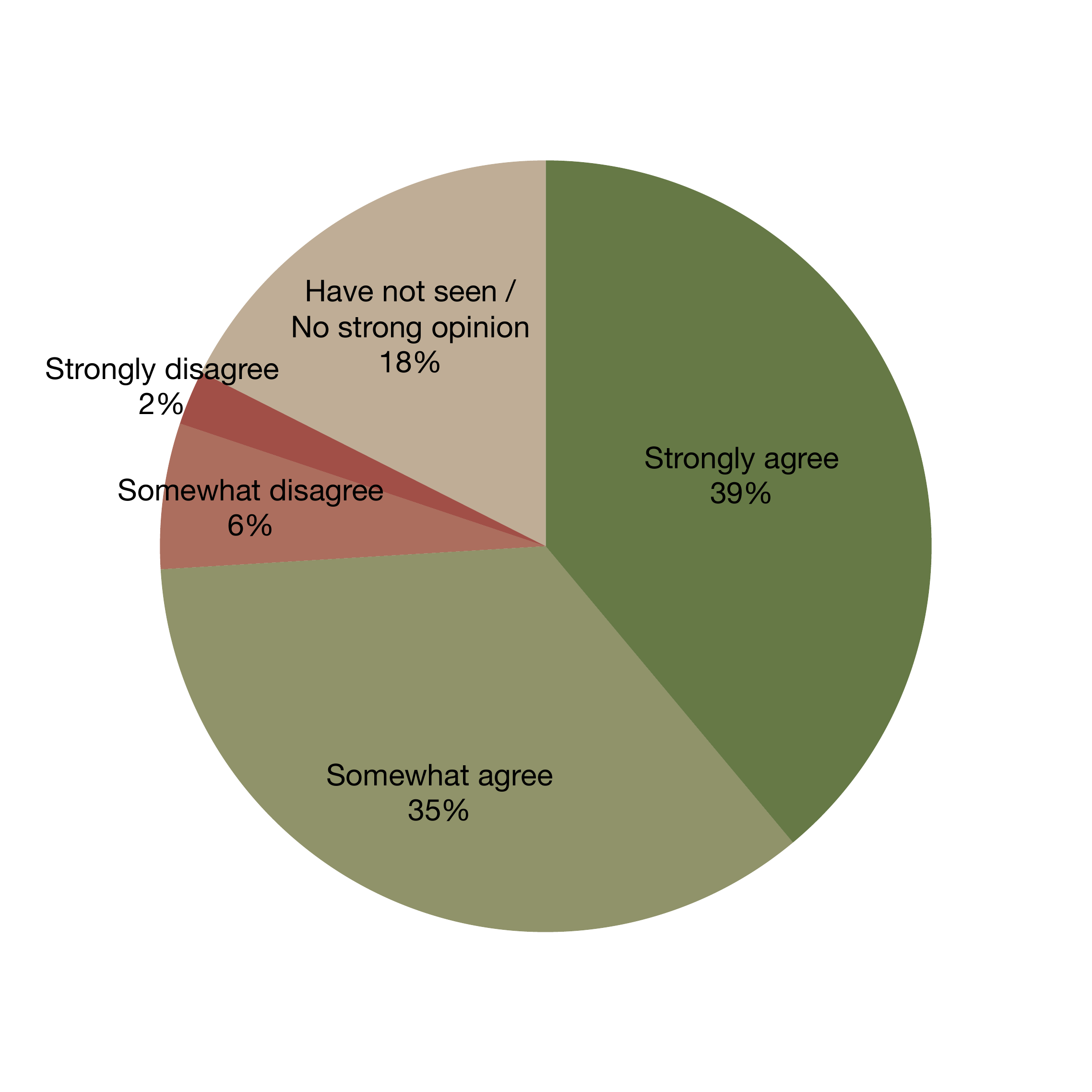}
    \vspace{-2.5\baselineskip}
    \caption{Responses to Question~\ref{q:declaration} ($n=262$).}
    \label{fig:declaration}
\end{figure}

\begin{question}\label{q:dept}
    Is your department doing enough to address AI? 
\end{question}

\begin{figure}[H]
    \centering
    \vspace{-\baselineskip}
    \includegraphics[width=\linewidth]{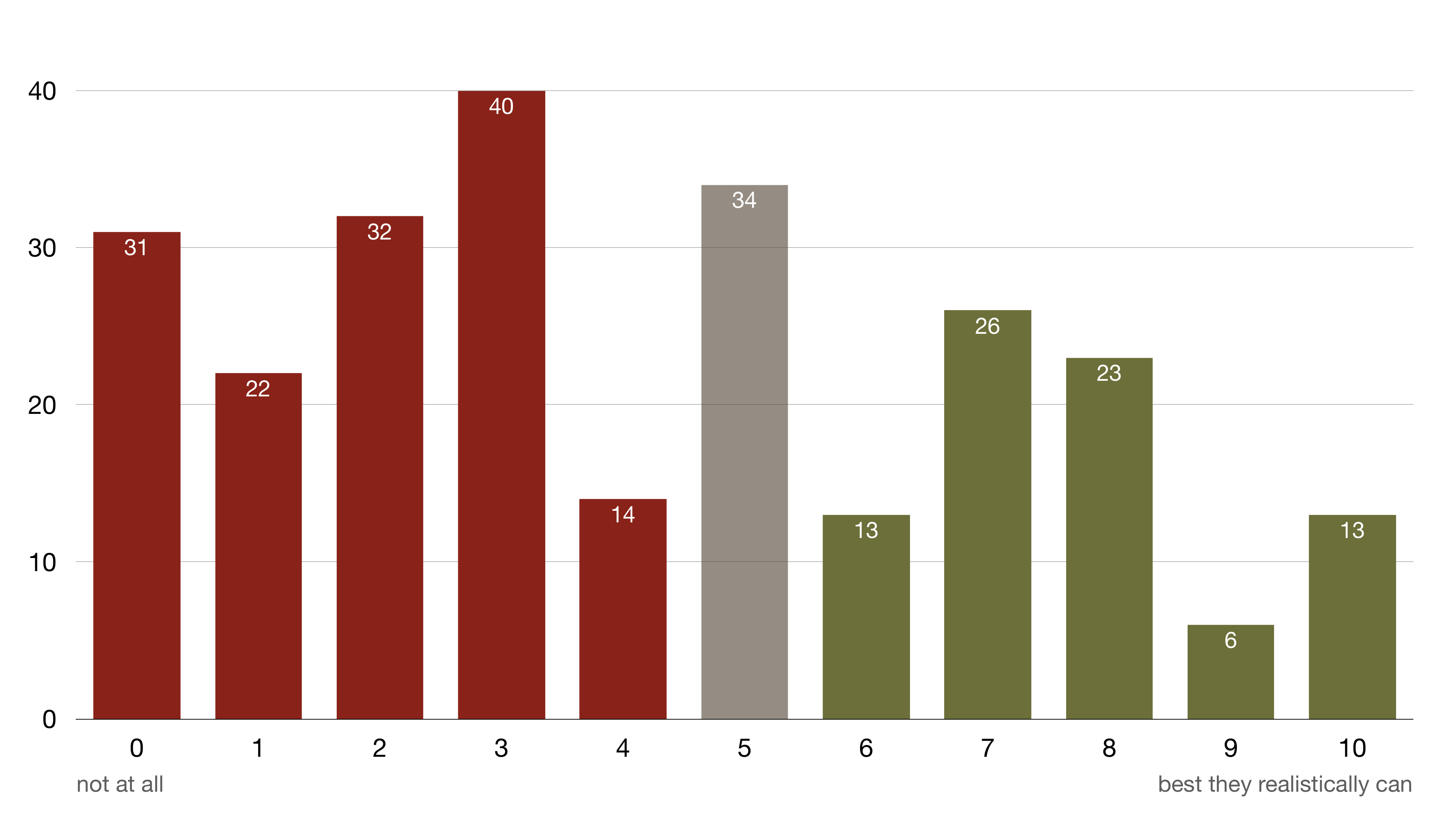}
    \vspace{-2\baselineskip}
    \caption{Responses to Question~\ref{q:dept} ($n=254$).}
    \label{fig:departments}
\end{figure}

\begin{question}\label{q:overall}
Which option comes closest to summarizing your general opinion on AI in math? 
\end{question}

\begin{figure}[H]
    \centering
    \vspace{-\baselineskip}
    \includegraphics[width=.85\linewidth]{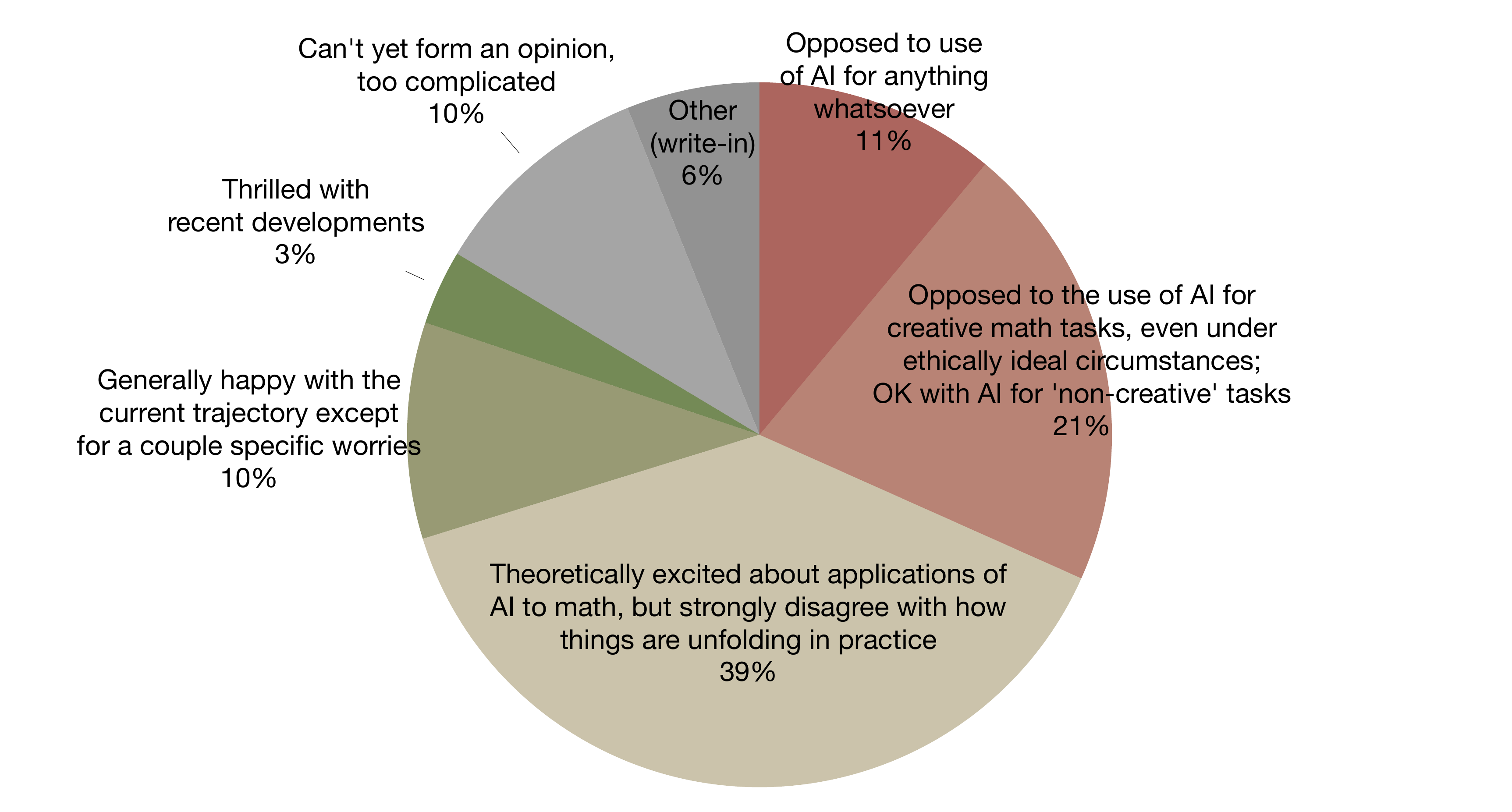}
    \vspace{-.5\baselineskip}
    \caption{Responses to Question~\ref{q:overall} ($n=262$).}
    \label{fig:opinion}
\end{figure}

\section{Conclusion}

Particularly striking to us was Question~\ref{q:continue}: `Concretely, has the probability of you continuing in academic math increased or decreased over the last couple of months?' 69.9\% of respondents said that the probability has decreased. What does this figure mean for the future of the field? 

As journal editors, hiring committees, faculty members, and leading figures in the field make decisions that affect the trajectory of our field, we urge them to take into account the opinions and worries of early career researchers. Math needs motivated, supported graduate students and postdocs. What can you do to motivate and support them better? 

We would love to hear from you with any thoughts you have on this survey. Feel free to email the authors: dorawood@mit.edu, abrande@mit.edu, hannamul@mit.edu. 

\printbibliography

\appendix

\bigskip
\bigskip
\begin{center}
    \textsc{Appendix}
\end{center}
\medskip

Our survey included several free response questions, where we invited respondents to elaborate on their viewpoints and experiences. Respondents were able to opt in or out of having their free responses anonymously quoted. 

To avoid making this document too long, we aren't recording every single free response; we chose a random sample of around half. A few of the quotes have been lightly edited for clarity/brevity/grammar. We've also lightly sorted responses by topic/subject within sections to help with readability.  

There are a lot of free responses below, and some of them are quite long. Nonetheless, we encourage you to take a look; we found many of the responses thought-provoking, surprising, and even moving. 

\section{Free Responses to `Any more thoughts on [your career trajectory and the job market?]'}

Funding and hiring pressures:
\begin{enumerate}
    \item `Ominous portents for funding.'
    \item `I worry about AI as a justification for hiring less tenure-track faculty and allocating less university/government funding to research in mathematics.'
    \item `If you're a university admin deciding between giving a professor $2,000$ worth of tokens or a $40,000$/year postdoc grant, one of those is a lot cheaper, and the companies who profit off of that decision insist that it's just as good. I have no faith in admin to make the right decision.'
    \item `With AI development rapidly becoming extremely proficient at research math, I worry that grant availability will decrease (as institutions begin to view human mathematicians as irrelevant in comparison).'
    \item `I think it will make [getting a job in academia] harder only if the number or tenure track jobs drastically decreases because of a change in society structure overall.'
\end{enumerate}

Uncertainty and unpreditability: 

\begin{enumerate}
    \item `It'd be great if institutes gave guidance instead of platitudes on this subject.'

    \item `It's hard to predict how the job market will change in academia. I find it equally plausible that the AI revolution increases the demand for mathematicians (there need to be more of us to keep up with AI) as that it decreases it (AI can do (most of) what we do, and as a result there will be less funding available).'

    \item `I'm not as worried about getting a tenure track role as I am about the future of the field in general. Say I get a TT role next year, great, but what about 5-10 years from now? what value will I bring in my TT role, if any?'

    \item `I think it’s extremely hard to predict what the job market will be like in several years. I could imagine scenarios where it gets easier or harder.'
    
    \item `I think it is nearly impossible to make any accurate prediction about what math will look like in the next few years.'

    \item `The uncertainty of the field is worrying for anyone at a critical career stage now, applying for postdoc or faculty positions.'

    \item `Given the speed of the development of AI, it's not clear at all to imagine what the academia becomes, how academia changes by the time I graduate.'

    \item `What will remain for us to do, what would jobs look like, research or otherwise, the uncertainty is real, and none of us can predict it; it almost feels like going with the flow.'
    
    \item `My worries about jobs in academia are somewhat synonymous with my worries about society more broadly (which is more important).'
    
\end{enumerate}

Institutional norms and leadership:
\begin{enumerate}
    \item `I am very concerned that the standard expectation for a job will be dozens of papers now, cranked out by people that don't fully understand them, rather than slower more thoughtful work.'

    \item `Math used to have 2 sides, solving and explaining. It seems like the human role is shifting to be just explaining. So specializing more in teaching and pedagogy seems increasingly important and may buy better job security.'

    \item `For me, the biggest source of stress in this regard is the uncertainty about the evaluation criteria. If they remain unchanged, I am afraid that the academic job market will turn into a big class struggle between those who can affort the most powerful models and those who can't.'

    \item `Too many things are in flux at the moment. We are lacking high-level leadership from across the academic community: from those who administer graduate programs, from those on various committees (graduate admissions, hiring, tenure), from those on editorial boards at journals, and from those who oversee public funding processes. Right now, individuals voice their own opinions in public, but we need to hear from those who can speak officially on behalf of institutions involved in the research process. (Right now, these include departments, universities, journals, and funding agencies.)'

    \item `I think new paradigm must be built, before that the academic career in AI era seems very unsure.'

    \item 
`Despite me not changing my plans in academia, I see a lot of decline in young researcher's willingness to work in academia. I think that as a community we need to protect young researchers. Here are some concrete things:
1) Even if ChatGPT can do this, but it can be a good project, give it to a Ph.D. student. My first ever paper happened because an older graduate student needed someone to do simple but tedious asymptotic computations in a generalized case of his previous paper. It would not be a hard task for AI, but I got to have a research project, learn a lot, and moreover had an older graduate student mentor that helped me throughout my Ph.D. Small projects matter for students to learn things. 
2) I think a lot of people already talked about valuing exposition and clearness. 
3) I would advise journals to not publish an article if an author cannot submit a video lecture explaining the paper. (In my opinion would also make the review process better). I think it is definitely very important that even if a proof is AI generated, someone in the community understands it, and actually can integrate it into the research area in the future. 
4) I think we should as a community should not accept proofs generated by AI companies without them writing the paper better and citing the relevant techniques. They can even pay someone to edit the paper, they can hire a whole team if needed - that creates jobs and then someone actually reads it. But we as a community should not be burdened with AI slop.'
\end{enumerate}

Meaningfulness of the job and value of human mathematicians:

\begin{enumerate}
    \item `The leading problem for me is that I don't enjoy doing the ``maximum AI" version of mathematics so much as the activities I am used to, so I would may prefer to leave if this is the case. This also leads to a small decrease in hireability, since, even if I \emph{could} compete with my peers using AI tools to a strong extent, I would rather not.'
    \item `It is difficult to do research and go to seminars when it feels like none of the individual work will matter much in the future.'
    \item `Math will become a game for AI agents. It’s the dawn for a new age, when human efforts become unimportant.'
    \item `I find myself reaching for it to prove technical lemmas more than I would like. I worry that it is likely capable of solving entire research problems, but I refrain from asking it to do so, preferring to guide the proof myself and using it to fill in details. I also worry that all of the papers that I've written pre-AI are cheapened: I spent countless hours on them, but many can probably be ``one-shotted" at this point.'
    \item `We will need to adapt, but it may be difficult to convince people why we need mathematicians after this.'
    \item `As someone in applied math, I think that a good understanding of interdisciplinary things and being able to apply such things in innovative ways will go a long way. I am quite comfortable with betting that AI won't replace that.'
\end{enumerate}

Personal career decisions and contingency plans: 
\begin{enumerate}
\item `At the present moment, the main career shift I am considering is moving into AI safety if I become sufficiently worried about existential AI risk and believe I can make a difference.'

\item `It hasn't changed my plans for now, but long term I think it's now more likely I'll be forced to leave math eventually due to AI. My plan is to ignore this and proceed with math for now, but be prepared to change careers if/when it comes to this down the line.'

\item `I'm concerned about what will happen with the jobs and if it is necessary to start studying something more applied.'

\end{enumerate}

\section{Free Responses to `Please share as much as you're willing about how AI has impacted your own resesarch.'}

Getting `scooped' by AI:
\begin{enumerate}
\item `In one incident, we hosted a visitor (a graduate student from another school) in our lab and presented some work in progress. Following this meeting, the visitor privately communicated to the PI an AI-generated proof of one of the problems. Within the lab, we had committed to working on the problem in the ``traditional" way, with little AI involvement (and certainly none for idea generation). While it was not surprising that the problem could be resolved by AI, it was unnerving for a visitor to prompt AI with the problem and show us the results without flagging that they would do so.'

\item `I have had a moderate result scooped by another group using AI (as usual it is difficult to tell how much AI was involved). This was rather unfortunate.
I also have a couple papers which were essentially write-ups of AI-discovered proofs. I think that the results in these papers are nice, and they are good for my CV, but I am not proud of them.'

\item `I guess as one more concrete suggestion for the community:
It is true that at least some people in the community would ask around before starting a problem to avoid competition and scooping. I think we as a community should do this even more often.'

    \item `Fortunately AI has not impacted my research yet, but I am concerned that it could in the future. I have become more reluctant to share my research (even pre-prints) on social media such as LinkedIn, and will probably not be sharing ongoing work as openly as previously. I have never fed any of my original research ideas into an AI, but only used it to brainstorm different broad/general strategies for proving technical conditions such as e.g., measurability.'

    \item `AI has negatively impacted my research because I am now afraid to share my ongoing work with people, which I love doing. Also, I have seen my collaborators have their papers killed by concurrent AI work, for instance following the nonsofic group announcement and the $\theta(p_c) = 0$ announcement. This makes the current environment very tense for me as a non-AI user, because I'm not sure what problems to work on except those that I come up with myself, which aren't formulated in the literature, yet I still fear that someone will just try to one-shot them. It somehow makes it difficult to think and commit to a project.'

\end{enumerate}

Impacts on personal motivation and growth:

\begin{enumerate}
    \item `From an individual perspective, I now feel coming up with technical proofs is meaningless. I once spent two weeks found a technical proof for some lemma, but it turns out to be incorrect and AI just found another much shorter proof immediately. This made me upset, and then decided to not do this by myself anymore. I turn to focus on more conceptual part now. Not sure if this is good or bad.'

    \item 'I don't spend as much time struggling to understand an argument in a paper. Sometimes I'm lazy and just ask AI to digest it for me and give me the step by step. I genuinely feel this is making me stupider, in the sense that every time I decide to skip the struggle of trying to understand just by reading and working things out, it gets harder next time I try to do it. To me this outweights the positive things, because there's no point on being a mathematician if you can't do mathematics.'

    \item `When I work on a problem from the very start, I need to analyze previous approaches people have used and then identify the real obstacles that stopped others to move forward. It’s precisely through these winding paths and dead ends that I build the correct intuition for the problem. This would eventually give me the ability to write down my own proof.

When AI gives me a logically rigorous proof directly, it’s like building a bridge between peaks of two mountains. I can walk across it and I'm sure that every plank is solid, but I’ve never touched the rock walls of the valley with my own hands, nor do I understand why the bridge had to be built right there. I know the proof is correct, but I don't know why it is correct.

Of course, even before AI, there is no single person that can re-derive every proofs on their own. But at the most time, we could find those who had walked the path themselves, and through human-to-human communication, the understanding is passing through the community. That, to me, is one of the main functions of our community.'

\end{enumerate}

Research incentive structures: 

\begin{enumerate}
    \item `I feel like the current atmosphere is driving academia toward a picture where people quickly produce a huge number of papers with AI, but the authors don't have enough time to fully understand what they’ve generated. Even if someone wants to slow down to deeply reflect on and polish their work, our environment likely won't give people that kind of space and freedom anymore.
    
    When I read an AI-generated proof, even if I get every single detail, it’s still hard for me to feel like the insight truly belongs to me. Without going through that artisan-style mathematical exploration to discover the answer myself, can I really gain a deep understanding? I don't think I can really grow to be a real mathematician by only reading proofs rather than writing down my own.(At the same time, the LLM itself has no memory—even if it gained a sort of understanding during its thinking process, that understanding is forgotten right after, it doesn't get the growth either.)'

    \item `Generally, it feels risky to work on publicly known conjectures - I could sink two months into a problem, only for the problem to get one-shotted by AI. This feeling has already deincentivized me from working on at least two projects. 

A specific instance: a fully AI-generated paper came out that was very relevant to one of my projects. It seems critical for me to understand the results of this paper if I want to make progress on my own project (eg, to make sure I'm not repeating their work). But the dense AI writing is very unclear (I can't even parse what the main theorem statement says) and NONE of the human authors are able to provide much clarity. As such, progress on my own project has been stalled for several weeks now.'

\item `Harder to justify incremental steps for the sake of machinery/theory building when the test problem at the end is bashed out.'

\item `One of the major game-changers for grad students is publishing, which can impact the way you think, increase your confidence, your public image, give you a chance to network, present your results, grow as a researcher. However, if everything related to publishing is now as good as a random coin flip, then where does the growth kick in? The temptation to not struggle on understanding a concept or working out a problem over the course of few days/weeks is real, and you cannot resist this temptation particularly because someone out there with weaker concepts and close to no understanding of the field can now produce the same results as you.'

\item `One aspect is advisor expects us to use it heavily, so taking the time to learn things the old fashioned way / move slower is less of an option. 
I also think it can make one lazy if not careful.
And it dramatically increases the odds of being scooped; certainly you can't put anything in a paper related to open problems or future directions anymore. It is very disheartening. In general I agree with points raised by others that the purpose of math should be more about building a general framework for understanding, not just rushing to prove theorems. And AI is not going in the right direction for this.'

\item `As a PhD student, AI is incredibly demotivating. I do not use AI as I believe this is a critical time in my development as a mathematician and off-loading thinking to AI at this stage will hurt my ability to solve problems. However, many other students use AI and have more publications than myself, which makes them perhaps more appealing for academic jobs. More importantly, it does not seem to me that people are being truthful with the extent of their AI use in disclosures.'

\end{enumerate}

Positive/mixed impacts of AI on research: 

\begin{enumerate}
    \item `Positive impact: I can understand papers much faster, because when I have a question, most times (when AI does not hallucinate), I don't have to wait until I meet with my advisor or any other expert to get it resolved. It is fairly good (even though there's a lot of room for improvement) for literature review, for example when wondering whether something has been proven before. It is also good for learning (as a non-native speaker) how do people typically state certain things.'

    \item `The day-to-day impact has been mostly positive; I don't get stuck nearly as often, because I have an endless stream of suggestions to keep me going. It's helped me feel less daunted by unsolved problems, because I have a clear starting point: getting ideas from AI. But I have missed a few obvious observations because they weren't in AI's suggestions. If you're given 5 ideas, it's hard to keep looking for a 6th instead of choosing among the 5.'

    \item `AI has increased my research productivity significantly, both in terms of speed and strength of results.'

    \item `AI has sped up my progress, and allowed me to attack different projects at the same time. I spend less time working on minute details/arguments and more thinking about big picture and how the smaller pieces fit together, as well as problem formulation. At the same time, I feel that the AI might soon overtake me even in those tasks in the next year or so. So paradoxically, even though I get more things done now, I am much more pessimistic about the future.'

    \item `In terms of helping my research, it is probably as helpful if not more helpful than my advisor. But still, sometimes it is hard to say if an idea is generated purely by AI or if it occurs to me through my discussion with AI (in an ``emerging" kind of way).'

\item `My research consists of writing a lot of code, and obviously AI is far better than me at this. It optimized my code dramatically (e.g., introducing parallel programming by batch). However, I feel it often over-engineers, leading to excessive useless code that I'm unable to review. For any task I give AI, there are many details that I don't think of at first, and LLM will do it in its own way, which deviates from what I think many times. AI has improved my code writing efficiency dozens or hundreds of times. I can write tests and verify an idea quickly. However, the feeling of losing control is depressing sometimes. That's why I feel both demotivated and motivated at the same time.'

\item `AI found errors in some of my writing that would have taken me much longer to spot without it. It has greatly sped up my ability to look for known results in the literature (e.g. textbooks which are 1000s of pages long can have key technical lemmas found by AI very quickly. This is extremely useful).'

\item `So far I have only used it for typo/error fixing in writing and some small improvements to technical lemmas, which has been mostly positive. I have mostly avoided having the problems I am interested in one-shotted but I don't expect that to last.'

\item `I’ve used AI to solve a few research questions I was interested in at this point, and I am genuinely excited to be able to see and understand solutions to these questions. However, I think really struggling with a problem is a very fun and rewarding aspect of doing math which leads to deep understanding, and I have felt depleted of that recently.'

\item `LLMs in mathematics are rapidly accelerating the pace of discovery and output. It is no longer possible to publish partial results on a problem, even if the problem in question is original and broadly interesting. There is too great a chance that someone may quickly publish a paper resolving the full problem, preventing publication. As we are, at least historically, judged by our publication records, this discourages the open nature of mathematical discourse.

On the other hand, the minutia of arguments are now much less of a focus during research. Where the proof of a lemma may have taken a few days in the past, it can now be done in an hour. We are free to think about the aesthetics of our work more, including how to organize and streamline ideas. As we generalize our ideas and push them farther and farther, the time it takes to rewrite material suitable to this generalization is orders of magnitude smaller.'
\end{enumerate}

Other thoughts:
\begin{enumerate}
    
    \item `I've given it a shot. I've toyed around with it, quizzed it on some topics that I know, and given it a toy lemma to try and prove for one of my projects; a lemma of the genre that I had been previously told it would be good at doing. It asserted that it was correct (which it was), asserted that it was well known (which I don't think it was), asserted there were citations in the literature (which it failed to produce), and when I asked it for a proof, the proof that it gave was faulty. It was, in my opinion, a complete and utter waste of time and money.

However, I think that AI and its reputation will have an outsized influence on my research career anyway. I'm an early career mathematician and have not made a name for myself; people will make the assumptions that they make, and unfortunately, the assumptions that people are making these days seems to be that everyone is using AI unless demonstrated otherwise, and I am afraid of that reputational harm.'

\item `I often find that AI tools slow me down, while using language that implies it is speeding me up. It will make me think there's a single line of code that can do what I want, but then it admits later that it used nonexistent syntax. It will tell me inaccurate definitions, sticking by them even when probed further, until I specifically say that it is wrong. I have given it code that was calculating something inaccurately before, and while it can spot the error in my syntax, it completely misses the inaccuracy in what my code is trying to compute. I have come to have very little confidence in anything that an AI model generates, and I often find that it takes me more time to verify its accuracy than to just figure things out the semi-old-fashioned way (referencing professors, books, and online forums). Not to mention, the latter option is significantly more fun.'

\item `AI only really entered my radar around July of this year, and I have really only been wrapping up my project and haven't started anything new, so I have been trying to avoid AI usage in the completion of this project. Sometimes when looking things up on Google, Gemini will ``accidentally" give me something useful, and I also ask ChatGPT to summarize ideas, conduct literature search, etc. Once I start a new project, maybe I will be a little more bold in asking LLMs for ideas/proofs of technical lemmas/etc, because it feels like I have to in order to keep up.'

\end{enumerate}

\section{Free Responses to `Any more thoughts on [AI affecting your research]?'}

\begin{enumerate}
    \item `It's been quite stressful to focus on work properly. Labs using my dream problems to train their kraken is certainly not pleasant. On the other hand, it forced me to reflect a bit deeper about my values and such.' 
\item `AI should never be used to write a paper. AI slop papers are going to do possibly irreparable damage to mathematics. Part of the reason that I have not been that discouraged so far is that I work in a low profile field where few people have taken to churning out slop. If/when that changes, my motivation may well plummet.'
\item `A problem I have run into is wishing for my pre-prints to be Open Access, but not wishing to consent to AI-models training on my work. I would ideally like for the different arXiv licenses or the community standards for using CC-BY licenses to be reconsidered to better protect authors legally from having their work stolen by AI, while still keeping research open access to the community.'

\item `Sometimes in mathematical research (if you're lucky), you can get the sense of making progress on a big mathematical problem within a few years. You might not know the exact details of how this will flesh out, but you sense that some of your ideas will make progress down the line and perhaps solve the problem altogether.

However, these days there is a good chance that if you embark on such a journey, you will get ``scooped" by some AI-model. I am anxious about how this will affect my motivation to pursue long-term projects. The idea of having an AI-model solving a problem that I have been working on for a long time is devastating. This has not yet really happened in my own field of mathematics, but I very much feel for the people for which it has.'

\item `I’ve just wrapped up a few projects that were done mostly without AI (the main ideas were mine or a collaborators). I haven’t started using AI for research very heavily, and I’m a bit worried that math won’t be as much fun for me in the future if I become reliant on AI. I’d like to say I won’t use it much in future projects for this reason, but I’m not entirely sure.'

\item `I feel more sad about the overall landscape of research given AI, but also more pressure to be spending time working and producing lots of results just to hopefully stay in math now.'

\item `Multiple times, from different professors, I was presented with an AI generated proof to check. Although it was certainly not their intention, I found that quite insulting, because it implicitly downgrades my role from collaborator to machine checker.
I feel AI is magnifying the differences (in this case of expertise) and inequalities within academia. If we keep thinking in productivity terms, I don't see how this revolution could be enjoyable, at least for those of us who do not have a permanent position.'

\item `While the academy may have strong ethical views regarding AI use and how it impacts the community, we also need to be mindful, especially on public platforms, that we remain accountable to the expectations of our funding sources, which indirectly includes the general public (via taxpayer funds). This includes responsible consumption of funds, which includes appropriate use of tools to conduct research efficiently. A decision for or against AI use must acknowledge the views from underlying funding sources. (Whether we respect those views completely is a subsequent matter.)'

\item `I know many people feel negatively about the role of AI in mathematics. While I agree we need to critically consider how these developments affect our work, and how to best move forward, I believe much of this resistance is motivated by emotions instead of rationality (and I say this despite having read the rational justifications for the negative outlooks, which I am not invalidating). Many mathematicians are grieving and adjusting, and I believe many people will process this grief and feel more positive as time goes by.'

\item `In general I think it's pretty cool that I can ask an AI system to try and implement some idea of mine to solve some problem. Usually the exchange is like this:
1) It tells me that it tried but either found a counterexample or gives me a ``target" lemma to prove which is almost harder than the thing I wanted to prove. 
2) I reassess my approach, based on some insight I gain from its response, modify it, and repeat. 

Sometimes it leads to genuinely useful ideas/progress and sometimes it makes clear to me that I either know nothing about the nature of the problem or that I'm really stubborn. I know that it's possible to simply run these models with a harness autonomously and, after many hours of compute, solve something like the Komlos conjecture, but I'm not too interested in doing that right now.'

\item `I would strongly prefer to never use AI in math, but over the past few months I have felt a competitive pressure: some mathematicians are using AI to speed up research, and if I don't want to be left behind I feel pressure to do this too. I think if we had a field-wide agreement about norms for using AI, I would strictly abide by the norms and perhaps use AI even less than the norms would allow, but since we don't have that there is a lot more uncertainty about how much/little AI use is necessary to succeed in research, and personally I don't want to use it any more than strictly necessary due to its adverse effects on the environment, society, and math itself.'

\item `I think I'm relatively relaxed about AI taking over this part of my mathematical activity, because personally pushing the frontier forward was, while exciting, never the main point. I just wanted to know, and however I get to know is fine.'

\end{enumerate}

\section{Free Responses to `Feel free to elaborate or push back on the previous question [about which concerns raised by the community worry you most]'}

The importance of AI safety/ethics beyond math: 

\begin{enumerate}

\item `The implications on math are an extremely small part of the overall AI issue right now. Seems like a large fraction of informed and thoughtful people believe we're on track to all literally die in the next few years, which would obviously be a much bigger deal than the rest of the things on this list. Hence why I feel weird checking any except the last of these boxes (although I do think that many of these would be concerns if we didn't have much bigger issues).'

\item `I think concerns about what these AI developments mean for math specifically are parochial. This is a tide sweeping over all human activity and the leading order term in determining our future, both short and long-term. We can sort out the situation with math once it has become clear whether we will die, survive as mankind, or as something else entirely.'

\item `Many of the bullet points would probably worry me much more if not for the AI safety concern. I think that one trumps all in the sense it is hard to say I'm ``worried whether math will be rewarding in ten years" when I'm more primitively worried whether there will be people to do math.'

\end{enumerate}

General competence and understanding of mathematicians: 

\begin{enumerate}
    \item `I just want to say that for me the most important as a community is ``A longterm decrease in human mathematical expertise/deep understanding of mathematics" and I think this will be a huge crisis if things continue the way they seem to be going. Also doing math will get way lonelier.'
    \item `One additional worry I have is the following: I think that the current propaganda of ``AI is good" will result in a cohort of students believing that they should use AI when they absolutely shouldn't, and as a result, we will have a generation of widespread incompetence. This will cause massive, lasting damage to mathematics even after AI is addressed.'
\end{enumerate}

Pushback / reasons for optimism: 
\begin{enumerate}
    \item `Everyone seems to agree that publishing new results will no longer be valued as much going forward, as it may just be the work of AI. I am relieved by this, as it will kill the toxic publish-or-perish culture in academia, which was the source of much ``slop" (yes, it existed before AI). I actually think humans will become more mindful about the things they spend their time on, and also believe quality of teaching and exposition will be more highly valued and therefore improve.'
    \item `I think it’s unlikely that no one will want to be a mathematician or even that there will be a significant brain drain from math in the scenario where all concrete problems are solved by AI, I think it will just be a different crowd. I also think the social aspects of math might get amplified, since the problem solving will take the backseat.'
    \item `Even before AI, academics were talking about the negative impact on researchers and their work of academic institutions pushing for quantity over quality. Similarly, students have been cheating on their work for a long time, looking things up rather than asking their teachers. There are wider social issues that encourage both of these. I often find myself drawing (imperfect) comparisons between AI and calculators in math. Calculators never stopped people who wanted to learn math from taking the time to understand problems that calculators can answer quickly. I do think there will always be people who want to understand.'
\end{enumerate}
Other thoughts: 
\begin{enumerate}
    \item `This has been the most interesting question from the survey. My belief now is that it will not be possible for the world to ``lose" all its mathematicians in the next year - even if AI develops strictly superhuman mathematical capabilities, we mathematicians will still report in to work, and hire other mathematicians, at least for a little while. My primary concern is that it will no longer be a ``fun" job that we get hired for - that there will be no humans (at least myself) capable of proving anything that an AI cannot. And if our job is reduced to digestion (or communication etc.) I don't think that will be very pleasant. Even worse, how will we be able to hire? I understand that there are some proposals there, but frankly as an early career mathematician I have not yet seen any that I trust. (Of course this is layered already onto a frankly depressing job market situation).'

    \item `Please do not use the word ``digestion". You are pandering to AI companies by using it, since it's their x accounts introduced this lexicon. We have not used this word to describe the process of understanding the proof.'

\item `I fear that the conversation about how our work will shift towards ``proof digestion" is mostly led by people who already have tenure, and are asking themselves what they will do once they're outclassed in terms of theorem proving. I have a hard time imagining that anyone will be willing to fund proof digestion on a large scale, and it seems we're deluding ourselves slightly in thinking this is what being a mathematician will become.'

\item `I think the lack of 'entry-level' problems could become redundant in a couple years time if we lack any problems humans can do that AI can't.'
\end{enumerate}

\section{Free Responses to `Any more thoughts on [the longterm health of the field]?'}
\begin{enumerate}
    \item `The breakneck pace of change is exhausting, and these worries are genuinely affecting my mental health.'
    \item `I don't even think I could make the argument to a highschool student as to why they should pursue a maths degree. Perhaps the only good argument is that I can't think of a different degree that might be more AI proof.'
    \item `The way I frame it is that certain AI uptake regimes will deteriorate our working conditions. We are not going to be able to attract people to research mathematics if we cannot offer good working conditions - time and space to think, a way of valuing research work that is value-pluralistic rather than solely about originality, and yes, credit and attribution for one's unique contributions. Right now, the atmosphere is extremely tense and uncertain, and the rhetoric is sometimes explicitly hostile to mathematicians who don't embrace AI. This isn't a good thing.'
    \item `I think that there could be a positive effect on teaching standards. I hope that if research skills become less valuable, then more hiring decisions will be based on ability to explain difficult concepts in an engaging and clear way.'
    \item `These developments are being discussed too much in the abstract as a wave of progress happening \emph{to} us, rather than many concrete decisions that individuals and companies are taking, mostly companies whose financials are extremely murky and dubious.'
    \item `Because so many variables are in flux, there is a possibility that none of the practices we are familiar with today exist in any recognizable form in the medium-term future. We will need to be careful in deliberating about the future state of certain practices in the community (or, at least, we may need to do so at much more fundamental and abstract levels). It is possible that the practices which we can describe and debate today will be supplanted altogether.'
    
\end{enumerate}

\section{Free Responses to `Final thoughts on the general issue of AI in math?'}

Broader ethics of using AI:

\begin{enumerate}
    \item `At present AI is used in so many unethical ways (deepfakes and CSAM, plagiarism, academic cheating, environmental issues, dehumanization of humanity, racism and other implict biases, ...) that it itself has become unethical technology which I cannot justify using in any capacity, regardless of its genuine usefulness. I'd be willing to use AI if all or at least the vast majority of its ethical problems are solved by some sort of regulation. But for now, metaphorically speaking, every time I make use of AI, an innocent person gets shot.'

    \item `The hypothetical ethically ideal circumstances do not exist. The businesses and people behind these tools do not have the best interest of mathematics, education more broadly, the Earth, or the next generation of humanity in mind. On top of this, LLMs are trained on human research, walked through questions by human researchers, motivated by humans, corrected by humans, with the accuracy of what they say determined by humans. I see no reason AI accomplishments shouldn't be credited to humans. But there are people with financial motive to present their tool as an all-powerful, separate-from-human thing.'

    \item `My frustrations with AI are my frustrations with current arguably monopolistic state of the US tech industry in general.'

    \item `I think that theoretically, AI could be a great tool for math if 1.\ AI use was limited to professional academics/people who could use it responsibly, in particular not minors (who won't learn how to think) or potential bad actors, 2.\ AI was not privatized (which is connected to 1), 3.\ AIs did not guzzle so much water (which would likely not be a significant issue if 1 and 2 were satisfied), and 4.\ there were uniform norms across our discipline for AI use.'

    \item `I have no moral qualm with the use of AI by ethical, trained mathematicians to better their research in our current world, but outside of the context of mathematics if I had a kill-switch for these tools I would press it in a heartbeat.'

    \item `Beyond mathematics, I think AI will have detrimental consequences for society. Even disregarding the extinction-case scenarios that have been discussed recently in the media, I worry that the decline in critical thinking skills arising from the widespread use of LLMs in  the general public, will be disastrous, and fundamentally incompatible with the idea of democracy (which relies on its citizens staying somewhat informed and being capable of forming independent opinions). Not to speak of general safety threats (cybersecurity, terrorism, biological weapons, etc.), the massive impact data centers have on the climate, and the massive concentration of power in tech-firms.'

    \item `I'm sympathetic to angst about what AI will do to the lives of math researchers, but kinda feel here like a dinosaur worrying about whether the incoming meteor's shadow is going to ruin the weather for my beach trip. It would be great if this could remain our biggest problem, but that's not really what the trends suggest...'

    \item `To your question of ``To what extent do you keep in mind broader, non-mathematical ethical issues in mind (eg: misinformation spread, military contracts) when you consider using AI for research?": I hadn't even thought to think about these things until reading this question. So now I will be thinking of these things more often.'
\end{enumerate}

Suggestions for the community (what are best case scenarios? silver linings? next steps?):

\begin{enumerate}
    \item `I think ideally AI should be used exclusively for refereeing via proof formalization and verification, and for educational purposes (literature review, understanding a topic/argument in a paper, etc.). 
I think this fever of productivity, for example publishing $>6$ paper per year, can only lead to a completely disconnected math community.'

\item `Thank you for this survey! I think it's essential that we coordinate a collective response to this major shift, both as students to the mathematical community and as future mathematicians to society.'

\item `I can't wait until there is a more universal agreement on the capabilities of math so that I don't have to shamefully only use AI by myself, and could use it directly in collaboration with my peers.'

\item `The math community needs to change significantly facing AI. We cannot simply say we're against AI or whatever, we need to design mechanism to help us coexist with AI. As for career decision, I'm waiting for this new mechanism to come out, and if it doesn't, I'll leave.'

\item `I've noticed professors often pushing AI on students more than they want it (and also referring them to AI when the student was asking for human mentoring/resources), and this is very alarming.
One thing only briefly touched on here: We need to address how citing AI works in math. I heard someone suggest there should be a standardized section included in papers about exactly if/how AI was used. I also think there should be particularly harsh consequences for putting out plagiarized work in any research that used AI to set a precedent of standards (and ideally have the risk discourage such AI use).'
\item `As communal stances from the academy emerge, we must demonstrate to the general public that our decisions represent a sensible application of our labour force. This may mean even accepting that the academy does not necessarily deserve to exist in the way that it does today. We must uphold the trust placed upon us by those who fund the community.'

\end{enumerate}

Academia's incentive structures, funding and hiring: 

\begin{enumerate}
    \item `I think AI has exposed some bad habits of the math community that were already there before.'

    \item `The way I frame it is that certain AI uptake regimes will deteriorate our working conditions. We are not going to be able to attract people to research mathematics if we cannot offer good working conditions - time and space to think, a way of valuing research work that is value-pluralistic rather than solely about originality, and yes, credit and attribution for one's unique contributions. Right now, the atmosphere is extremely tense and uncertain, and the rhetoric is sometimes explicitly hostile to mathematicians who don't embrace AI. This isn't a good thing.'

    \item `There is a lot of pressure to be the first one to solve research problems. If AI takes a major role in math research, we should really try to understand the results, but as things stand there is incentive to try and get more results instead.'

    \item `I'm optimistic we can realign our incentives, but we need to do it now, because LLMs break our old incentive system. The whole economy of ``theorem credits" does not make sense anymore.'

    \item `I hope we can restructure around collective understanding and human exposition. My biggest worry is how maths will remain meritocratic in the future.'

    \item `I'm mostly concerned about the future of funding and hiring after AI (I’m also concerned with larger social impacts like the public losing critical thinking or human extinction, but I’ll focus on math issues). Most of my anxiety around AI comes from not knowing if it will make it impossible to find a job.
    
    As far as the use of AI in math, I feel most mathematicians would be happy with an AI that is at the level of a PhD student or coauthor who can help them when they get stuck or do busy work. However they aren’t ready for what to do when there’s no use even thinking about a problem since an AI will be 100x more successful. In some ways it is exciting that AI will let us learn many new things, but disappointing that there’s many problems I’ll never have the chance to be the one to solve. I don’t think this problem is unique to mathematicians. Many people will soon have to face the reality that an AI can do every part of their job better than them, and as a society we need to figure out what that means.
    
    I don’t think there’s any reasonable future where AI progress is immediately halted or AI is forbidden in mathematics. Even if AI stayed at its current capabilities, it would be extremely disruptive as more people get used to it. And if progress continues at the current rate, do any of my papers have value in a few years?
    
    I also want to mention that many senior researchers don’t realize that they’re in a privileged position because they’ve already proved their expertise before AI. Nobody’s going to fire a professor of 30 years with tenure. But how will future hiring committees know that I’m smart if an AI can replicate my papers? Why should a future university hire new postdocs or professors?'
    
    \item `My greatest fear for the future of mathematics as a discipline has much less to do with the replacement of human cleverness and ingenuity with large language models (or whatever future technology replaces it). I believe, as an intellectual pursuit among a community, the value of mathematics and its training will remain in perpetuity, even as its practice and training should evolve. 
    
    What concerns me most is the economic and administrative uncertainty that the rapid evolution of LLMs has brought, and no doubt will continue to bring. While it is very well and good to discuss the importance of human understanding of mathematics, and it is important, we must realize that the future of mathematics as a discipline may ultimately not be up to us. It is the accountants, administrators, boards of directors, and various other apparatchiks who will ultimately decide whether mathematics receives funding, is required to train students, or receives office space. If, as I believe, mathematical expertise will be of growing importance in the research landscape and it is at the same time misunderstood, perverted, or otherwise devalued by the bean counters of university administration and funding organizations, mathematics as a discipline could only persist as the lowest of creatures, a roach trying not to be stamped out. 
    
    It is the responsibility of mathematicians then, not only to adopt the tools that will make ourselves the most useful, including but not limited to LLMs, but also to pursue avenues to make ourselves more directly and visibly useful to the broader scientific community.'
    
\end{enumerate}

More thoughts on the meaning of doing math:
\begin{enumerate}
    \item `I think most mathematicians could have been a lot of other things. Most of those things would probably pay more and provide greater prestige. I think most of them chose mathematics because they loved doing math. But what does it mean to do math now? I think doing math without AI is fundamentally different from doing math with AI.'

\item `At the end of the day, I want to understand math. AI can only assist me in this. It is possible, or maybe even likely, that the presence of AI makes it harder to get a job whose sole goal is understanding math, but I probably won't end up with such a job anyways, independent of AI existing. It's probably an unpopular opinion influenced by my relative mediocrity as a mathematician, but I find the thought of quitting math research due to AI being able to do it better than me logically inconsistent with the idea that math has inherent beauty. Those who want to quit because they think math has been ``ruined" seem to imply that math being elegant or cool or interesting is conditional, and I disagree with that.'

\item `I think it's fundamentally wrong to talk about `The' purpose of mathematics. As humans we don't only care about fundamentally understanding things, only as mathematicians is that our goal. Also AI will clearly continue progressing faster than humans. A good analogy to make there is how there is so much content on youtube and feed algorithms are so good now, that if you talk to 20 people you might find barely any overlap between what they consume. Exaggerating somewhat, it used to be that there are 4 tv channels and everyone watched at least 2 of them. I think math interests will similarly progress to be incredibly niche like this. There will be far more results that are known to even smaller groups of people. I suppose that's good though if somehow global population continues to increase without us destroying the planet.'

\item `In an ideal world (putting eg job security, government funding for math, and potentially catastrophic environmental degradation aside), I would be extremely excited about AI applications to mathematics: it is clear that AI has the potential to supercharge our collective understanding of mathematics, and frankly I love learning new mathematics.

Ironically, in the long term I expect AI to decrease our collective understanding of mathematics, even if it becomes a mathematical oracle that can answer any question we have. This is because human understanding will become devalued; fewer people will choose to study mathematics, because I expect opportunities for funded study and research in mathematics are going to decrease drastically.'

\end{enumerate}
Other thoughts: 
\begin{enumerate}
    \item `Good luck.'
    \item `AI made me confused in lots of things.'
    \item `I think a large missing piece from this survey was even asking what we thought about the competence of AI at all. I've given it a shot and I find it deeply incompetent at original research; it seems capable only of plagiarism. My personal hypothesis is that the real mathematicians using AI in their own work are not really getting out of it what they say they are, but it's just a tendency that some mathematicians have to downplay their own contributions finally finding a convenient scapegoat to give credit to. That is, beyond my many ethical hangups over AI; I also just don't think that it's actually remotely useful.'

    \item `My department is flailing. Behind closed doors many of them are angry. But they have not committed to anything. They don’t know what to do.'

    \item `Clearly a useful tool to learn if used carefully, but using it to crank out papers in a PhD seems like a good way to learn nothing.'
        
    \item `Thank you for doing this! If there are any practical ways that I can help, please let me know!'

    \item `The declaration ‘A Severe Misalignment of AI in Mathematics?’ came out stupid and sad.  It was at that moment that many people in the AI community lost respect for field medal and pure mathematics as a whole.  Navier Stokes was solved.  Instead of staying up all night reading, digesting, relishing the proof content, those mathematicians focused on writing that Ludditean piece of joke.'
    
    \item `I see it as an inevitability. I feel that it would be better if it never existed, but since it does exist, there's no going back. It will be used ethically and unethically, and we can't stop it. Those who don't use it will fall behind, so everyone will have to adapt. Such is life.'

    \item `I'm getting a bit tired of tenured professors pontificating about the value of mathematical understanding. A lot of these declarations are fluff and are detached from reality. They also serve absolutely no purpose in alleviating the anxiety of the next generation.'
\end{enumerate}

\end{document}